\documentclass[aos,MSNbibl,dvips]{arximspdf}
\usepackage{mathbh}
\usepackage{graphicx}

\doi{10.1214/12-AOS1015} 
\volume{40}
\issue{3}
\pubyear{2012}
\firstpage{1578}
\lastpage{1608}

\makeatletter
\newcommand{\rrvert}{\vert}
\newcommand{\llvert}{\vert}
\newcommand{\eqref}[1]{(\ref{#1})}
\newcommand
\prob{\mathbb{P}}
\def\argmax{\operatorname{argmax}}
\renewcommand{\hat}{\widehat}
\newcommand\indicator{\mathbh{1}}
\renewcommand{\d}{\mathrm{d}}

\newcommand{\E}{\mathbb{E}}
\newcommand{\N}{\mathbb{N}}
\renewcommand{\P}{\mathbb{P}}
\newcommand{\R}{\mathbb{R}}
\newcommand{\holder}{\sigma}

\newtheorem{theo}{Theorem}[section]
\newtheorem{cor}{Corollary}[section]
\newtheorem{lem}{Lemma}[section]
\newproclaim{remark}{Remark}[section]
\makeatother

\begin{document}
\begin{frontmatter}

\title{The limit distribution of the $L_\infty$-error of Grenander-type
estimators}
\runtitle{The $L_\infty$-error of Grenander-type estimators}

\begin{aug}
\author[a]{\fnms{C\'ecile} \snm{Durot}\ead[label=e1]{cecile.durot@gmail.com}},
\author[b]{\fnms{Vladimir N.} \snm{Kulikov}\ead[label=e2]{vladimir.kulikov@asr.nl}}
\and
\author[c]{\fnms{Hendrik P.} \snm{Lopuha\"a}\corref{}\ead[label=e3]{h.p.lopuhaa@tudelft.nl}}

\runauthor{C.~Durot, V. N. Kulikov and H. P. Lopuha\"a}

\affiliation{University of Nanterre, ASR Nederland and Delft
University of Technology}

\address[a]{C. Durot\\
Universit\'e Paris Ouest Nanterre La D\'efense\\
200, avenue de la R\'epublique\\
92001 Nanterre Cedex\\
France\\
\printead{e1}}

\address[b]{V. N.~Kulikov\\
ASR Nederland\\
Archimedeslaan 10\\
3584 BA Utrecht\\
The Netherlands\\
\printead{e2}}

\address[c]{H.~P. Lopuha\"a\\
DIAM\\
Delft University of Technology\\
Mekelweg 4\\
2628 CD Delft\\
The Netherlands\\
\printead{e3}}
\end{aug}

\received{\smonth{11} \syear{2011}}
\revised{\smonth{5} \syear{2012}}

%
\begin{abstract}
Let $f$ be a nonincreasing function defined on $[0,1]$. Under standard
regularity conditions, we derive the asymptotic distribution of the
supremum norm of the difference between $f$ and its Grenander-type
estimator on
sub-intervals of $[0,1]$.
The rate of convergence is found to be of order $(n/\log n)^{-1/3}$ and
the limiting distribution to be Gumbel.
\end{abstract}

%
\begin{keyword}[class=AMS]
\kwd[Primary ]{62E20}
\kwd{62G20}
\kwd[; secondary ]{62G05}
\kwd{62G07}.
\end{keyword}

\begin{keyword}
\kwd{Supremum distance}
\kwd{extremal limit theorem}
\kwd{least concave majorant}
\kwd{monotone density}
\kwd{monotone regression}
\kwd{monotone failure rate}.
\end{keyword}

\end{frontmatter}

\section{Introduction}
\label{secintroduction}
After the derivation of the nonparametric maximum likelihood estimator
(NPMLE) of a monotone density and
a monotone failure rate by Grenander~\cite{grenander1956}, and the
least squares estimator
of a monotone regression function by Brunk~\cite{brunk1958}, it has
taken some time before the distribution
theory for such estimators entered the literature.
The limiting distribution of the NPMLE of a decreasing density on
$[0,\infty)$ at a fixed point in the interior of the support,
has been established by Prakasa Rao~\cite{prakasarao1969}.
Similar results were obtained for the NPMLE of a monotone failure rate
in~\cite{prakasarao1970}
and for an estimator of a monotone regression function in~\cite{brunk1970}.
Woodroofe and Sun~\cite{woodroofe-sun1993} showed that the NPMLE of a
decreasing density is inconsistent at zero.
The behavior at the boundary has been further investigated in~\cite
{kulikov-lopuhaa2006,balabdaouietal}.
Smooth estimation has been studied in~\cite{mammen1991}, for monotone
regression curves, and
in~\cite{vdvaart-vdlaan2003} for monotone densities; see also~\cite
{eggermont-lariccia2000}
and~\cite{anevski-hossjer2006}.
The limit distribution of the NPMLE of a decreasing function in the
Gaussian white noise model was
obtained in~\cite{zilberburg2007}.
Related likelihood ratio based techniques have been investigated
in~\cite{banerjeewellner2001,pal2009}.

Groeneboom~\cite{groeneboom1985} reproved the result in~\cite
{prakasarao1969} by introducing a new approach based on inverses.
This approach has become a cornerstone in deriving pointwise
asymptotics of several shape constrained nonparametric estimators,
for example, for the distribution function of interval censored
observations (see~\cite{groeneboom-wellner1992})
or for estimators of a monotone density and a monotone hazard under
random censoring (see~\cite{huang-wellner1995});
see also~\cite{huang-zhang1994} for the limiting distribution of the
NPMLE of a monotone density under random censoring
and~\cite{lopuhaanane2011} for similar results on isotonic estimators
for a monotone baseline hazard in Cox proportional hazards model.
The limit distribution of these estimators
involves an argmax process $\{\zeta(c)\dvtx c\in\mathbb{R}\}$ connected
with two-sided Brownian motion with a parabolic drift.
This process has been studied extensively in~\cite{groeneboom1989},
where it is also claimed that the approach
based on inverses should be sufficiently general to deal with global
measures of deviation, such as the $L_1$-distance or the supremum
distance between the estimator and the monotone function of interest.
Indeed, the limiting distribution of the $L_1$-distance between a
decreasing density and its NPMLE was obtained in~\cite
{groeneboom-hooghiemstra-lopuhaa1999},
and a similar result can be found in~\cite{durot2002} in the monotone
regression setup.
These results were extended to general $L_k$-distances in~\cite
{kulikov-lopuhaa2005} and~\cite{durot2007}.
In~\cite{durot2007}, the limiting distribution of $L_k$-distances is
obtained in a very general framework
that includes, among others, the monotone density case, monotone
regression and monotone failure rate.

Little to nothing is known about the behavior of the supremum distance.
In~\cite{jonker-vdvaart2001}, the rate of the supremum distance is
established in a semi-parametric model for censored observations,
and it is suggested that the same rate should hold in the monotone
density case.
In~\cite{hooghiemstra-lopuhaa1998} an extremal limit theorem has been
obtained for suprema of the process $\zeta(c)$ over increasing intervals.
However, a long-standing open problem remains, although this problem
has important statistical applications: what is the limiting
distribution of the supremum distance
between a monotone function and its isotonic estimator?
Indeed, while pointwise
confidence intervals for a decreasing density, a~monotone regression function or a monotone hazard are available using
the limiting distribution of the isotonic estimator at the fixed point,
nonparametric confidence bands have remained a formidable challenge;
they could be built if the limiting distribution of the supremum distance
between a~monotone function and its isotonic estimator were known.
It is the purpose of this paper to settle this question in the same
general framework as considered in~\cite{durot2007}.
The precise construction of a nonparametric confidence band requires
additional technicalities that are
beyond the scope of the present paper.
It is only briefly discussed here, and details are deferred to a
separate paper.

We consider Grenander type estimators $\widehat{f}_n$ for decreasing
functions $f$ with compact support, say $[0,1]$.
These are estimators that are defined as the left-hand slope of the
least concave majorant of an estimator for the primitive of $f$.
This setup includes Grenander's~\cite{grenander1956} estimator of a
monotone density,
Brunk's~\cite{brunk1958} estimator for a monotone regression function,
as well as the estimator
for a monotone failure rate under random censoring, considered in~\cite
{huang-wellner1995}.
We obtain the rate of convergence for the supremum of~$|\widehat
{f}_n-f|$ over subintervals of $[0,1]$.
The rate is shown to be of the order $(n/\log n)^{-1/3}$, even on
subintervals that grow toward~$[0,1]$,
as long as one stays away sufficiently far from the boundaries,
so that the inconsistency at the boundaries (see, e.g.,~\cite
{woodroofe-sun1993}) is not going to dominate the supremum.
The rate that we obtain coincides with the one suggested in~\cite
{jonker-vdvaart2001}
for Grenander's~\cite{grenander1956} estimator for a decreasing density,
but it is now proven rigorously in a more general setting under optimal
conditions on the boundaries of the intervals
over which $\sup|\widehat{f}_n-f|$ is taken.
Moreover, we show that the rate $(n/\log n)^{-1/3}$ is sharp.
Our main result is Theorem~\ref{theolawfn}, in which we show that a
suitably standardized supremum of $|\widehat{f}_n-f|$
converges in distribution to a standard Gumbel random variable.

Our results are obtained following the same sort of approach as that used
in~\cite{groeneboom1985,huang-wellner1995,groeneboom-hooghiemstra-lopuhaa1999,durot2002,durot2007}, among others.
We first establish corresponding results for the supremum of the
inverses of~$\widehat{f}_n$ and~$f$,
and then transfer them to the supremum of $\widehat{f}_n$ and $f$ themselves.
A major difference with deriving asymptotics of $L_k$-distances is,
that in these cases one can benefit from the linearity of the integral
and handle several approximations
pointwise with Markov's inequality.
This is no longer possible with suprema.
With suprema, to transfer results for inverses to results for $\hat f_n$,
a key ingredient is a precise uniform bound on the spacings between
consecutive jump points of $\hat f_n$.

The paper is organized as follows.
In Section~\ref{secmainresults}, we list the assumptions under which
our results can be obtained
and state our main results concerning the rate of convergence and the
limiting distribution of $\sup|\widehat{f}_n-f|$.
We also briefly discuss the construction of confidence bands.
We formulate corresponding results for the supremum distance between
the inverses of~$\widehat{f}_n$ and~$f$
in Section~\ref{secinverse}.
This is the heart of the proof, which is carried out in Section~\ref
{secproofsUn}.
Finally, in Section~\ref{secproofratelawfn}, we
provide a uniform bound on the spacings between consecutive jump points
of~$\hat f_n$ and then transfer the results obtained in
Section~\ref{secinverse}
for the inverses of~$\widehat{f}_n$ and~$f$ to the supremum distance
between the functions themselves.

To limit the length of the paper, the rigorous proofs of several
preliminary results needed
for the proofs in Sections~\ref{secproofsUn} and~\ref{secproofratelawfn}
have been put in a~supplement~\cite{durotkulikovlopuhaa2012}.

\section{Assumptions and main results}
\label{secmainresults}
Based on $n\ge2$ independent observations, we aim at estimating a
function $f\dvtx [0,1]\to\R$ subject to the constraint that it is
nonincreasing. Assume we have at hand a cadlag (right continuous with
finite left-hand limits at every point) stepwise estimator $F_n$ of
\[
F(t)=\int_0^t f(u)\,\d u,\qquad  t\in[0,1],
\]
with finitely many jump points.
In the case of i.i.d.~observations with a~common density function $f$,
a typical example is the empirical distribution function with $n$
discontinuity points located at the observations.
In the following, we shall consider the monotone estimator $\hat f_n$
of~$f$ as defined in~\cite{durot2007}, that is,~the estimator $\hat
f_n$ is the left-hand slope of the least concave majorant of~$F_n$ with
\[
\hat f_n(0)=\lim_{t\downarrow0}\hat f_n(t).
\]
As detailed in Section \ref{secunifrate} below, this definition
generalizes well-known monotone estimators, such as the Grenander
estimator of a nonincreasing density, or the least-squares estimator of
a monotone regression function.
It should be noted that $\hat f_n$ is nonincreasing, left-continuous
and piecewise constant.
We are interested in the limiting behavior of the supremum distance
between the monotone estimator and the function $f$.

\subsection{Uniform rate of convergence}
\label{secunifrate}
We first show that the rate of convergence of $\hat f_n$ to $f$ in
terms of the supremum distance is of order $(\log n/n)^{1/3}$.
To this end, we make the following assumptions.
Unless stated otherwise, for a function~$h$ defined on $[0,1]$,
we write $\Vert h \Vert_\infty= \sup_{t \in[0,1]} \vert h(t) \vert$.
\begin{longlist}[(A3)]
\item[(A1)]
The function $f$ is decreasing and differentiable on $[0,1]$ with
\[
\inf_{t \in[0,1]} \bigl|f'(t)\bigr| > 0\quad \mbox{and}\quad \sup_{t \in[0,1]}
\bigl|f'(t)\bigr| < \infty.
\]
\item[(A2)] Let $B_n$ be either a Brownian bridge or a Brownian motion.
There exist $q\ge4$, $C_q>0$, $L\dvtx [0,1]\to\R$ and versions of $F_n$ and
$B_n$ such that
\[
\P \bigl(n^{1-1/q} \bigl\Vert F_n -F - n^{-1/2}B_n
\circ L \bigr\Vert_\infty >x \bigr) \le C_qx^{-q}
\]
for all $x\in(0,n]$. Moreover, $L$ is increasing and differentiable on
$[0,1]$ with $\inf_{t \in[0,1]}L'(t)>0$ and $\sup_{t \in[0,1]}
L'(t)<\infty$.
\item[(A3)] There exists $C_0>0$ such that for all $x>0$ and $t=0,1$,
\[
\E \Bigl[\sup_{u\in[0,1], x/2\le|t-u|\le
x}\bigl(F_n(u)-F(u)-F_n(t)+F(t)
\bigr)^2 \Bigr] \le\frac{C_0x}{n}.
\]
\end{longlist}
These conditions are similar to the ones used in~\cite{durot2007}.
Assumption~(A1) is completely the same as the one in~\cite{durot2007}.
Assumption~(A2) is similar to (A4) in~\cite{durot2007}, but now we only
require $q\geq4$
and bounds on the first derivative of~$L$.
Here we can relax the condition on $q$, because in the current
situation the error terms
have to be of smaller order than $(n/\log n)^{1/3}$ instead of
$n^{1/2}$ in~\cite{durot2007}.
The existence of $L''$, as imposed in~(A4) in~\cite{durot2007}, is not
needed to establish Theorem~\ref{theoratefn}.
Finally, assumption~(A3) is equal to (A2$'$) in~\cite{durot2007}.
Assumption~(A2) in~\cite{durot2007} is no longer needed,
since we are able to obtain sufficient bounds on particular tail
probabilities with our current
assumptions~(A1)--(A2).
See Lemma~6.4 and also the proof of Lemma~6.10 in~\cite{durotkulikovlopuhaa2012}.

A typical example that falls into the above framework
is the problem of estimating a nonincreasing density $f$ on $[0,1]$.
Assume we observe i.i.d. random variables $X_{1}, X_2,\ldots,X_{n}$ with
common nonincreasing density function $f\dvtx [0,1]\to\R$, and let $F_n$ be
the corresponding empirical distribution\vadjust{\goodbreak} function. In this case, the
monotone estimator $\hat f_{n}$ of $f$ coincides with the Grenander estimator.
Assumption (A1) is equal to the ones in~\cite
{groeneboom-hooghiemstra-lopuhaa1999,kulikov-lopuhaa2005,durot2007},
and is standard when studying $L_k$-distances between~$\widehat{f}_n$ and~$f$.
The existence of a~second derivative of $f$ is not needed to obtain
Theorem~\ref{theoratefn}.
In the monotone density model, assumption (A2) is satisfied for all $q>0$,
with $L=F$ being the distribution function corresponding to $f$
and $B_n$ a~Brownian bridge,
due to the Hungarian embedding of~\cite{komlosmajortusnady1975}.
From Theorem 6 in~\cite{durot2007} it follows that assumption (A3)
holds in the monotone density model.
Another example that falls into the above framework
is the problem of estimating a monotone regression function.
Assume for instance that we observe $y_{i}=f(i/n)+\varepsilon_{i}$,
$i=1,2,\ldots,n$, where the $\varepsilon_{i}$'s are i.i.d. centered
random variables with a finite variance $\sigma^2$, and $f\dvtx [0,1]\to\R$
is nonincreasing. Let $F_{n}$ be the partial sum process given by
\[
F_{n}(t)=\frac1n\sum_{i=1}^ny_{i}
\indicator_{i\leq nt}.
\]
In this case, the monotone estimator $\hat f_{n}$ of $f$ coincides with
the Brunk estimator.
Assumption (A1) is equal to the ones in~\cite{durot2002,durot2007} and
is standard when studying $L_k$-distances in this model. Assumption
(A2) is satisfied for all $q\geq2$ such that $\E|\varepsilon_{i}|^q<\infty$ with $L(t)=\sigma^2t$ and $B_{n}$ a Brownian motion,
due to embedding of~\cite{sak1985}. Thus, (A2) is satisfied in the
above regression model provided $\E|\varepsilon_{i}|^4<\infty$. From
Theorem 5 in~\cite{durot2007} it follows that assumption (A3) holds in
the above regression model.
Other examples of statistical models that fall in the above framework,
with corresponding~$q$ and~$L$, are discussed in~\cite{durot2007}.

The uniform rate of convergence of $\hat f_n$ to~$f$
for general Grenander-type estimators is given in the following theorem.
%
\begin{theo}\label{theoratefn}
Assume \textup{(A1), (A2)} and \textup{(A3)}.
Let $(\alpha_n)_n$ and~$(\beta_n)_n$ be sequences of positive numbers
such that
%
\begin{equation}
\label{eqalphanratefn} \alpha_n \ge K_1
n^{-1/3}(\log n)^{-2/3} \quad\mbox{and}\quad \beta_n\ge
K_2 n^{-1/3}(\log n)^{-2/3}
\end{equation}
for some $K_1,K_2>0$ that do not depend on $n$.
Then,
\[
\sup_{t\in(\alpha_n,1-\beta_n]}\bigl|\hat f_n(t)-f(t)\bigr|=O_p \biggl(
\frac{\log
n}{n} \biggr)^{1/3}.
\]
\end{theo}
The rate in Theorem~\ref{theoratefn} coincides with the one found
for the maximum likelihood estimator in a semi-parametric model for
censored data
by Jonker and van der Vaart~\cite{jonker-vdvaart2001},
who suggest that this rate should also hold for Grenander's~\cite
{grenander1956} estimator for a decreasing density.
They consider $\alpha_n\gg n^{-1/3}(\log n)^{1/3}$ and $\beta_n$ constant,
which is a slightly stronger assumption than the one in Theorem~\ref
{theoratefn}.
Note that condition~\eqref{eqalphanratefn} in Theorem~\ref{theoratefn} is sharp.
If $\alpha_n=n^{-\gamma}$, for some $1/3<\gamma<1$, then
$n^{(1-\gamma)/2}(\widehat{f}_n(\alpha_n)-f(\alpha_n))$ converges in
distribution,
according to Theorem~3.1(i) in~\cite{kulikov-lopuhaa2006},
so that
\[
(n/\log n)^{1/3}\bigl|\hat f_n(\alpha_n)-f(
\alpha_n)\bigr|\to\infty.
\]
In fact, for sequences $(\alpha_n)_n$ such that $n^{1/3}(\log
n)^{2/3}\alpha_n\to0$,
it can be shown similarly that $(n\alpha_n)^{1/2}\{\hat f_n(\alpha_n)-f(\alpha_n)\}$ converges in distribution,
which would yield
$(n/\log n)^{1/3}|\hat f_n(\alpha_n)-f(\alpha_n)|\to\infty$.

\subsection{Limiting distribution}
Whereas the previous theorem only provides a bound on the rate of convergence,
it is nevertheless crucial for deriving the actual asymptotics of the
supremum norm of $\hat f_n - f$ on suitable intervals.
For this purpose, we need an additional H\"older assumption on $f'$ and~$L''$.
\begin{longlist}[(A4)]
\item[(A4)]
The function $L$ in (A2) is twice differentiable and there exist
$C_0>0$ and $\holder\in(0,1]$ such that for all $t,u\in[0,1],$
%
\begin{equation}
\label{eqA4} \bigl|f'(u)-f'(t)\bigr|\le C_0|u-t|^\holder
\quad\mbox{and}\quad \bigl|L''(u)-L''(t)\bigr|\le
C_0|u-t|^\holder.
\end{equation}
\end{longlist}
The condition on $L''$ in assumption (A4) is a bit stronger than the
one in~\cite{durot2007}.
This is needed to guarantee that the difference between the values
of~$L''$ at~$t$
and its nearest point of jump of $\widehat{f}_n$ is negligible.
The condition on~$f'$ in assumption (A4) is the same as (4) in~\cite{durot2007},
who already observed that the existence of $f''$, as assumed in~\cite
{groeneboom-hooghiemstra-lopuhaa1999,kulikov-lopuhaa2005},
is no longer needed.
Note that in the monotone density model $L''=f'$, in which case (A4)
reduces to a H\"older condition on~$f'$ only.
In the monotone regression model, $L$ is linear so that~(A4) again
reduces to a H\"older condition on~$f'$ only.

In order to formulate the limit distribution, we need the following
definition:
%
\begin{equation}
\label{defzeta} \zeta(c) = \mathop{\argmax}_{t \in\R} \bigl\{ W(t+ c ) -
t^2 \bigr\} \qquad\mbox{for all }c\in\R,
\end{equation}
where $W$ is a standard two-sided Brownian motion on $\R$ originating
from zero, and argmax denotes the greatest location of the maximum.
For fixed $t\in(0,1)$, properly scaled versions of $n^{1/3}(\widehat
{f}_n(t)-f(t))$
converge in distribution to the random variable $\zeta(0)$
(see, e.g.,~\cite{prakasarao1969} or~\cite{groeneboom1985}).
Moreover, $\zeta$ serves as the limit process for properly scaled versions
of $n^{1/3}(\widehat{U}_n-g)$ (see, e.g., Theorem~3.2 in~\cite
{groeneboom-hooghiemstra-lopuhaa1999}),
where~$\hat U_n$ and $g$ are the inverse functions of $\hat f_n$ and
$f$ respectively, as defined in Section \ref{secinverse} below.
Properties of the process $\{\zeta(c),c\in\R\}$ can be found in~\cite
{groeneboom1989};
for example, the process $\{\zeta(c),c\in\R\}$ is a stationary process.
According to Corollary~3.4 in~\cite{groeneboom1989}, the tails of the
density~$\mu$ of $\zeta(0)$
satisfy the following expansion:
%
\begin{equation}
\label{eqexpansionmu} \mu(t)\sim2\lambda|t|\exp\bigl(-2|t|^3/3-
\kappa|t|\bigr)
\end{equation}
as $|t|\to\infty$, where $\kappa$ and $\lambda$ are positive
constants.\vadjust{\goodbreak}

We now present the main result of this paper.
It states that the limit distribution of the supremum distance between
$\hat f_n$ and $f$, if properly normalized, is Gumbel.
By $x_n\gg y_n$ we mean $x_n/y_n\to\infty$, as $n\to\infty$.
%
\begin{theo}\label{theolawfn}
Assume that \textup{(A1), (A2), (A3)} and \textup{(A4)} hold.
Consider $0\leq u<v\leq1$ fixed.
Then, for any sequence of real numbers $(\alpha_n)_n$ and~$(\beta_n)_n$
both satisfying
%
\begin{equation}\label{eqcondan}
\qquad\alpha_n\to0, \qquad\beta_n\to0\quad
\mbox{and}\quad 1-v+\beta_n \mbox{, }u+\alpha_n\gg
n^{-1/3}(\log n)^{-2/3},
\end{equation}
we have that for any $x\in\R$,
%
\[
 \P \biggl( \log n \biggl\{ \biggl(\frac{n}{\log n}
\biggr)^{1/3} \sup_{t\in(u+\alpha_n,v-\beta_n]} \frac{\vert\widehat f_n(t) - f(t)\vert}{|2f'(t)L'(t)|^{1/3}} - \mu_n
\biggr\} \le x \biggr)
 \to \exp \bigl\{-\mathrm{e}^{-x} \bigr\}
\]
as $n \to\infty$,
where
%
\begin{equation}
\label{defmun} \mu_n= 1 - \frac{\kappa}{2^{1/3}(\log n)^{2/3}} + \frac1{\log n} \biggl[
\frac13\log\log n + \log(\lambda C_{f,L}) \biggr],
\end{equation}
with
%
\[
C_{f,L} = 2\int_u^v
\biggl(\frac{\vert f'(t) \vert^2}{L'(t)} \biggr)^{1/3} \,\mathrm{d}t,
\]
and $\lambda$ and $\kappa$ taken from~\eqref{eqexpansionmu}.
\end{theo}
Note that from Theorem~\ref{theolawfn}, with $u=0$ and $v=1$,
it follows that for convenient~$\alpha_n$ and~$\beta_n$,
\[
\biggl(\frac{n}{\log n} \biggr)^{1/3} \sup_{t\in(\alpha_n,1-\beta_n]}
\frac{\vert\widehat f_n(t) - f(t)\vert}{|2f'(t)L'(t)|^{1/3}}=1+o_p(1).
\]
Since both $f'$ and $L'$ are bounded from above and bounded away from zero,
this proves that there are positive numbers $C_1,C_2$ that depend only
on $f'$ and $L'$ such that
\[
C_1+o_p(1)\leq \biggl(\frac{n}{\log n}
\biggr)^{1/3} \sup_{t\in(\alpha_n,1-\beta_n]}\bigl\vert\widehat f_n(t) - f(t)
\bigr\vert\leq C_2+o_p(1).
\]
This means that the rate in Theorem \ref{theoratefn} is sharp.

\subsection{Confidence bands}
\label{secdiscussion}
Our main motivation for proving Theorem~\ref{theolawfn} is to build
confidence bands for a monotone function $f$.
Indeed, this theorem ensures that for any $x\in\mathbb{R}$, with probability
tending to $\exp(-\mathrm{e}^{-x})$, we have
\[
\bigl\vert\widehat f_n(t) - f(t)\bigr\vert\leq \biggl(\frac{\log n}{n}
\biggr)^{1/3}\bigl|2f'(t)L'(t)\bigl|^{1/3}
\biggl\{\mu_n +\frac{x}{\log n} \biggr\},
\]
simultaneously for all $t\in(u+\alpha_n,v-\beta_n].$
Combining this with either plug-in estimators of $f'$ and $L'$
or
bootstrap methods would provide a confidence band for~$f$,\vadjust{\goodbreak}
at the price of additional technicalities.
Indeed, the use of plug-in estimators for the derivatives $f'$ and $L'$
may lead to inaccurate intervals for small sample sizes $n$,
so that bootstrap methods should be preferable.
But it is known that the standard bootstrap typically does not work for
Grenander-type estimators;
see~\cite{kos2008,banerjee-sen-woodroofe2010}.
Thus, we shall use a smoothed bootstrap, which will raise the question
of the choice of the smoothing parameter.
In view of all this, we believe that the precise construction of a
confidence band is beyond
the scope of the present paper and is deferred to a separate
paper.\looseness=1

Note that the conditions of Theorem~\ref{theolawfn} do not cover the
supremum distance over the whole interval $[0,1]$.
However, this is to be expected.
For instance, consider the monotone density model.
This model is one of the examples that is covered by our general setup
(see Section~\ref{secunifrate}) and it is well known
that the Grenander estimator $\widehat{f}_n$ in this model is
inconsistent at 0 and~1
(e.g., see~\cite{woodroofe-sun1993}).
Therefore, a distributional result can only be expected if the supremum
is taken over subintervals of $[0,1]$
that do not include~0 and~1.
Let us notice, however, that we can obtain a confidence band for $f$ on
any sub-interval $(u,v]$ with fixed $u,v\in(0,1)$ (by considering
$\alpha_n=\beta_n=0$),
and that the largest interval on which our result allows to build a
confidence band is $(\alpha_n,1-\beta_n]$, where $\alpha_n\gg
n^{-1/3}(\log n)^{-2/3}$ and similarly, $\beta_n\gg n^{-1/3}(\log n)^{-2/3}$.
In order to obtain a confidence band on the whole interval $[0,1]$,
we would have to slightly modify the Grenander-type estimator $\hat f_n$
in order to make it consistent near the boundaries.
For instance, we conjecture that, if we consider either the modified
estimator in~\cite{kulikov-lopuhaa2006}
or the penalized estimator in~\cite{woodroofe-sun1993} instead of $\hat f_n$,
then the limit distribution of the supremum distance between this
modified estimator and $f$
over the whole interval $[0,1]$ is the same as the limit distribution
of the supremum distance
between $\hat f_{n}$ and $f$ over the largest interval allowed in
Theorem~\ref{theolawfn}.
Thus, such modified estimators would provide a confidence band for $f$
over the whole interval $[0,1]$.
As mentioned above, the precise construction of confidence bands is
deferred to a separate paper,
and we will do similarly with the precise study of modified estimators
at the boundaries.

\section{The inverse process}
\label{secinverse}
To establish Theorems~\ref{theoratefn} and~\ref{theolawfn}, we
use the same approach as
in~\cite
{groeneboom1985,groeneboom-hooghiemstra-lopuhaa1999,durot2002,durot2007}.
We first obtain analogous results (i.e., rate of convergence and limit
distribution)
for the supremum between the inverses of $\widehat{f}_n$ and~$f$,
and then transfer them to the supremum between the functions~$\widehat
{f}_n$ and $f$ themselves.
Let $F_n^+$ be the upper version of $F_n$ defined as follows:
$F_n^+(0)=F_n(0)$ and for every $t\in(0,1]$,
\[
F_n^+(t)=\max \Bigl\{F_n(t),\lim_{u\uparrow t}F_n(u)
\Bigr\}.
\]
Let $\widehat{U}_n$ denote the (generalized) inverse of $\widehat{f}_{n}$,
defined for $a\in\R$ by
$\widehat{U}_{n}(a)=\sup\{t\in[0,1]\dvtx  \widehat{f}_{n}(t)\geq a\}$,
with the convention that the supremum of an empty set is zero.
This is illustrated in Figure~\ref{figUn} below.
From Figure~\ref{figUn}, it can be seen that the value $t=U_n(a)$
maximizes $F_n^+(t)-at$,
so that
%
\begin{equation}
\label{defUn} \hat U_n(a)=\mathop{\argmax}_{t\in[0,1]} \bigl
\{F_n^+(t)-at \bigr\}.
\end{equation}
The advantage of characterizing the inverse process $\hat U_n$ by (\ref
{defUn}),
is that in this way, it is more tractable than the estimator $\hat f_n$
itself, as being the argmax of a relatively simple process.
It is the purpose of this section to establish results analogous to
Theorems~\ref{theoratefn} and~\ref{theolawfn} for the inverse
process.\looseness=1

\begin{figure}

\includegraphics{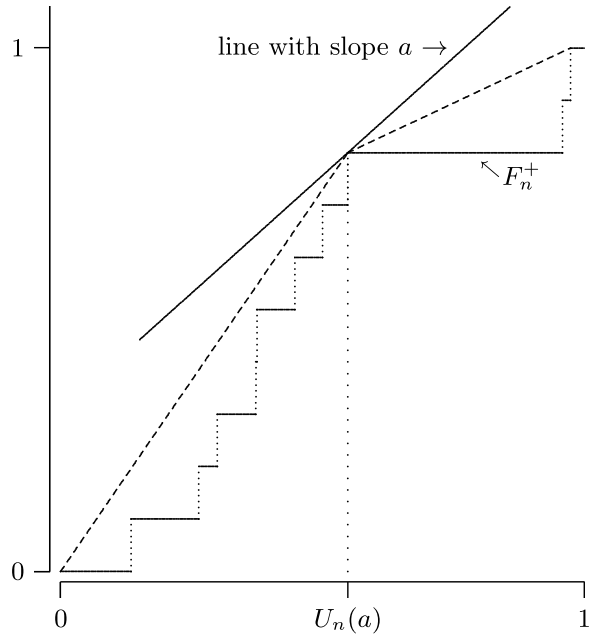}

\caption{The function $F_n^+$, its concave majorant (dashed) and a
line with slope $a$ (solid).}\label{figUn}\vspace*{3pt}
\end{figure}

Let $g$ denote the (generalized) inverse function of $f$.
In Theorems~\ref{theorateUn} and~\ref{theolawUn}, we give an
upper bound for the rate of convergence of~$\hat U_n$
to $g$, and an extremal limit result for the supremum distance
between~$\hat U_n$ and~$g$.
We derive the limit distribution of the supremum distance between~$\hat
U_n$ and~$g$ in Corollary~\ref{corlawUn}.
%
\begin{theo}\label{theorateUn}
Assume that \textup{(A1)} and \textup{(A2)} hold. Then
\[
\sup_{a\in\R} \bigl\vert\hat U_n(a)-g(a) \bigr\vert=O_p
\biggl(\frac{\log
n}{n} \biggr)^{1/3}.
\]
\end{theo}
%
\begin{theo}
\label{theolawUn}
Assume that \textup{(A1)}, \textup{(A2)} and \textup{(A4)} hold, and define for $a \in\R$ the
normalizing function
%
\begin{equation}
\label{defA} A(a)=\frac{|f'(g(a))|^{2/3}}{(4L'(g(a)))^{1/3}}.
\end{equation}
Let $0\leq u<v\leq1$ fixed, and let $(\alpha_n)_n$ and $(\beta_n)_n$
be sequences such that
$\alpha_n\to0$, $\beta_n\to0$ and $0\leq u+\alpha_n<v-\beta_n\leq1$
for $n$ sufficiently large.
Define
%
\begin{equation}
\label{defSn} S_n=n^{1/3}\sup_{a \in[f(v-\beta_n),f(u+\alpha_n)]} A(a)\bigl
\llvert \hat U_n(a) - g(a)\bigr\rrvert .
\end{equation}
Then
%
\begin{equation}
\label{eqLimitsupU} \prob ( S_n\le u_n ) \to \exp \biggl
\{ -2\tau\int_u^v \frac{|f'(t)|^{2/3}}{(4L'(t))^{1/3}}
\,\mathrm{d}t \biggr\}
\end{equation}
for any sequence $(u_n)_n$ such that $u_n\to\infty$ in such a way that
$n^{1/3}\mu(u_n)\to\tau>0$,
where $\mu$ denotes the density of $\zeta(0)$, as defined in~\eqref
{defzeta}.
\end{theo}
The expansion in~\eqref{eqexpansionmu} allows us to provide a
precise expansion of $u_n$
[see~(\ref{eqexpansionun})] and to derive the following corollary
from Theorem~\ref{theolawUn}. According to this corollary, the limit
distribution of $S_n$ is Gumbel.
%
\begin{cor}
\label{corlawUn}
Assume that \textup{(A1)}, \textup{(A2)} and \textup{(A4)} hold.
Let $S_n$ be defined by~\eqref{defSn}, with $0\leq u<v\leq1$, and
$\alpha_n,\beta_n$ satisfying the conditions of Theorem~\ref{theolawUn}.
Then, for all $x\in\R$,
\[
\P \biggl\{ \log n \biggl\{ \biggl( \frac{2}{\log n} \biggr)^{1/3}S_n-
\mu_n \biggr\} \le x \biggr\} \to \exp \bigl\{-\mathrm{e}^{-x}
\bigr\},
\]
where $\mu_n$ is defined by \eqref{defmun}.
\end{cor}
In order to transfer the results for $\widehat{U}_n$ to $\widehat
{f}_n$, we establish Lemma~\ref{lemslope-inverse}.
This lemma does require conditions on sequences $s_n=u+\alpha_n$ and
$t_n=1-v+\beta_n$
that are stronger than the ones in Theorem~\ref{theolawfn}.
However, once we have established the limit distribution for such sequences,
we will show that Theorem~\ref{theolawfn} can be extended to more
general sequences satisfying~\eqref{eqcondan}.

\section{\texorpdfstring{Proofs of Theorems~\protect\ref{theorateUn} and~\protect\ref
{theolawUn} and Corollary~\protect\ref{corlawUn}}
{Proofs of Theorems 3.1 and 3.2 and Corollary 3.1}}
\label{secproofsUn}
We suppose in the sequel that assumptions (A1) and (A2) are fulfilled,
and we denote by~$C$, $C_1$, $C_2,\ldots$ positive real numbers that
depend only on $q$, $C_q$, $f$, $L$ [and possibly on~$\holder$
under the additional assumption (A4)].
These real numbers may change from one line to the other.
We write $x\vee y=\max(x,y)$ and $x\wedge y=\min(x,y)$, for any real
numbers~$x$ and~$y$.

In order to deal simultaneously with the cases where~$B_n$ is a Brownian
bridge or a Brownian\vadjust{\goodbreak} motion [see assumption (A2)], we shall make use of
the representation
%
\begin{equation}
\label{eqrepresentation} B_n(t)=W_n(t)-\xi_n
t,\qquad t \in[0,1],
\end{equation}
where $W_n$ is a standard Brownian motion, $\xi_n\equiv0$ if $B_n$ is
a Brownian motion and $\xi_n \equiv W_n(1),$ a standard Gaussian
variable that is independent of $B_n$, in case $B_n$ is a Brownian bridge.
To prove Theorem~\ref{theorateUn}, we need some preliminary results
on the tail probabilities of $\widehat{U}_n-g$
and its supremum.
These results can be found in Supplement B in~\cite{durotkulikovlopuhaa2012}.
A first result, which is similar to Lemmas 2, 3 and 4 in~\cite
{durot2007}, is that
there exist $C_1 > 0$ and $C_2 > 0$ such that for all $a\in\R$ and $x>0$,
%
\begin{equation}
\label{1boundUn} \P \bigl( n^{1/3} \bigl\vert\hat U_n(a)-g(a)
\bigr\vert>x \bigr)\le \frac{C_1n^{1-q/3}}{x^{2q}}+2\exp\bigl(-C_2x^3
\bigr).
\end{equation}
In particular, for all $a\in\R$, this implies that
$\hat U_n(a)-g(a)=O_p(n^{-1/3})$.
See Lemma~6.4 in~\cite{durotkulikovlopuhaa2012}.
This is not sufficient to obtain Theorem~\ref{theorateUn},
but it will be used for its proof.

\begin{pf*}{Proof of Theorem~\protect\ref{theorateUn}}
Recall that $g(a)=1$ for all $a\le f(1)$, $g(a) = 0$ for $a \ge f(0)$
and~$\hat U_n$
is nonincreasing and takes values in $[0,1]$. Hence, we can write
%
\begin{equation}
\label{eqsuplef1} \sup_{a\le f(1)}\bigl|\hat U_n(a)-g(a)\bigr|= \bigl|\hat
U_n\bigl(f(1)\bigr)-g\bigl(f(1)\bigr)\bigr|
\end{equation}
and
%
\begin{equation}
\label{eqsupgef0} \sup_{a\ge f(0)}\bigl|\hat U_n(a)-g(a)\bigr|= \bigl|\hat
U_n\bigl(f(0)\bigr)-g\bigl(f(0)\bigr)\bigr|.
\end{equation}
This means that
\[
\sup_{a\in\R}\bigl|\hat U_n(a)-g(a)\bigr|= \sup_{a\in[f(1),f(0)]}\bigl|\hat
U_n(a)-g(a)\bigr|.
\]
Therefore, to prove Theorem \ref{theorateUn} it suffices to show that
\[
\sup_{a\in[f(1),f(0)]}\bigl|\hat U_n(a)-g(a)\bigr|=O_p \biggl(
\frac{\log n}{n} \biggr)^{1/3}.
\]
According to Lemma~6.5 in~\cite
{durotkulikovlopuhaa2012}, the bound in~\eqref{1boundUn} can be extended
such that for any $x>0$,
\begin{eqnarray*}
&&\P \biggl( \sup_{a\in[f(1),f(0)]}\bigl|\hat U_n(a)-g(a)\bigr| > x \biggl(
\frac{\log n}{n} \biggr)^{1/3} \biggr)
\\
&&\qquad \le \widetilde{C}_3 n^{1/3} \biggl(\frac{C_1n^{1-q/3}}{x^{2q}(\log n)^{2q/3}}+2n^{-C_2x^3}
\biggr),
\end{eqnarray*}
where $\widetilde{C}_3=C_3\{f(1)-f(0)\}$.
The latter upper bound tends to zero as $n\to\infty$ for all
$x>(3C_2)^{-1/3}$ since $q\ge4$ by assumption.
This completes the proof of Theorem \ref{theorateUn}.\vadjust{\goodbreak}
\end{pf*}

We suppose in the sequel that in addition to (A1) and (A2),
assumption~(A4) is fulfilled.
The first step in proving Theorem \ref{theolawUn} is to approximate
an adequately
normalized version of~$\hat U_n(a)$ by the location of the maximum of a
Brownian motion with parabolic drift.
To this end define
%
\begin{equation}
\label{defVn} V_n(a)= n^{1/3} \bigl( L\bigl(\hat
U_n\bigl(a^\xi\bigr)\bigr) - L\bigl(g(a)\bigr) \bigr),
\end{equation}
where
%
\begin{equation}
\label{defaxi} a^\xi = a - n^{-1/2} \xi_n
L'\bigl(g(a)\bigr)\qquad \mbox{for all $a\in\R$,}
\end{equation}
with $\xi_n$ taken from representation~\eqref{eqrepresentation}.
Then for $0\leq u<v\leq1$ and $\alpha_n,\beta_n$ satisfying the
conditions of Theorem~\ref{theolawUn},
we obtain
\begin{eqnarray*}
S_n\vee O_p(1) &=& \sup_{a \in[f(v-\beta_n),f(u+\alpha_n)]}
\frac{A(a) }{L'(g(a))}\bigl|V_n(a)\bigr| \vee O_p(1)
\\
&&{}+ O_p \bigl(n^{-\holder/2}(\log n)^{2/3}
\bigr)+O_p \bigl(n^{-1/6} \bigr),
\end{eqnarray*}
where $S_n$ is defined by~\eqref{defSn}, and $\sigma\in(0,1]$ is
taken from (A4).
See Lemma~6.6 in~\cite{durotkulikovlopuhaa2012}.

Next, we proceed with localization.
The purpose of this is that localized versions of~$V_n(a)$ and $V_n(b)$,
can be approximated by independent random variables, if~$a$ and~$b$
are in disjoint intervals that are suitably separated.
First note that the location of the maximum of a process is invariant
under addition of constants or multiplication by $n^{2/3}$.
Therefore, from~(\ref{defUn}) it follows that for all $a\in\R$ we have
%
\begin{equation}
\label{eqVn} V_n(a) = \mathop{\argmax}_{t \in I_n(a)} \bigl
\{W_{g(a)}(t) + D_n(a,t) + R_n(a,t) \bigr\},
\end{equation}
where
%
\begin{equation}
\label{defIna} I_n(a) = \bigl[n^{1/3} \bigl(L(0) - L
\bigl(g(a)\bigr)\bigr), n^{1/3} \bigl(L(1) - L\bigl(g(a)\bigr)\bigr)
\bigr]
\end{equation}
for every $s\in[0,1]$ fixed, $W_{s}$ is the standard Brownian motion
defined by
%
\begin{equation}
\label{defWga} W_{s}(t) = n^{1/6} \bigl\{ W_n
\bigl(L(s) + n^{-1/3} t \bigr) - W_n\bigl(L(s)\bigr) \bigr\}
\qquad\mbox{for $t\in\R$,}
\end{equation}
with $W_n$ defined by (\ref{eqrepresentation}), and
%
\begin{eqnarray}
\label{defRn} D_n(a,t) &=& n^{2/3} \bigl(F \circ
L^{-1} - a L^{-1} \bigr) \bigl(L\bigl(g(a)
\bigr)+n^{-1/3} t \bigr) \nonumber\\
&&{}- n^{2/3} \bigl(F\bigl(g(a)\bigr) - a g(a)\bigr),
\nonumber
\\[-8pt]
\\[-8pt]
\nonumber
R_n(a,t) &=& n^{2/3} \bigl(a - a^\xi \bigr)
\bigl(L^{-1}\bigl(L\bigl(g(a)\bigr) + n^{-1/3} t\bigr) - g(a)
\bigr)
\\
&&{} - n^{-1/6} \xi_n t + \tilde{R}_n(a,t) ,\nonumber
\end{eqnarray}
where $\xi_n$ is taken from representation~\eqref{eqrepresentation},
and for all $a$ and $t$,
%
\begin{equation}
\label{eqtildeRn} \bigl\vert\tilde{R}_n(a,t) \bigr\vert\le n^{2/3}
\bigl\Vert F_n - F - n^{-1/2} B_n \circ L
\bigr\Vert_\infty.
\end{equation}
For all $a\in\R$, we define the localized version of $V_n(a)$ by
%
\begin{equation}
\label{eqtildeVn} \tilde{V}_n(a) = \mathop{\argmax}_{t \in I_n(a):   \vert t \vert \le\log n} \bigl
\{W_{g(a)}(t) + D_n(a,t) + R_n(a,t) \bigr\}.
\end{equation}
We find that
\begin{eqnarray*}
&& \sup_{a \in[f(v-\beta_n),f(u+\alpha_n)]} \frac{A(a)
}{L'(g(a))}\bigl|V_n(a)\bigr|
\\
&&\qquad = \sup_{a \in[f(v-\beta_n),f(u+\alpha_n)]} \frac{A (b(a) )}{
L'(g (b(a) ))} \bigl\vert\tilde{V}_n(a) \bigr\vert
+ o_p(\log n)^{-2/3}
\end{eqnarray*}
for any $b(a)\in\R$ that satisfies $|a-b(a)|\le n^{-1/3}(\log n)^2$.
See Lemma~6.7 in~\cite{durotkulikovlopuhaa2012}.

Finally,
using the fact that, roughly speaking,
\[
D_{n}(a,t)\approx-\frac{\vert f'(g(a)) \vert}{2
(L'(g(a)))^2}t^2\approx-
\frac{\vert f'(g(b)) \vert}{2 (L'(g(b)))^2}t^2
\]
for all $b$ close enough to $a$, we bound $ \vert\tilde{V}_n(a) \vert$
from above and below by the absolute value of the
following quantities:
%
\begin{equation}
\label{defVntilde+} \tilde{V}_n^+(a, b) =\mathop{\argmax}_{t\in I_n(a): \vert t \vert\le\log n}
\biggl\{ W_{g(a)}(t) - \biggl(\frac{\vert f'(g(b)) \vert}{2 (L'(g(b)))^2} - 2\epsilon_n
\biggr) t^2 \biggr\}
\end{equation}
and
%
\begin{equation}
\label{defVntilde-} \tilde{V}_n^-(a, b) = \mathop{\argmax}_{t\in I_n(a): \vert t \vert\le\log n}
\biggl\{ W_{g(a)}(t) - \biggl(\frac{\vert f'(g(b)) \vert}{2 (L'(g(b)))^2} + 2\epsilon_n
\biggr) t^2 \biggr\},
\end{equation}
where $I_n(a)$ and $W_{g(a)}$ are defined in~\eqref{defIna}
and~\eqref{defWga}, $b$ is chosen sufficiently close to $a$,
and where $(\epsilon_n)_n$ is a sequence of positive numbers that
converges to zero as $n\to\infty$,
which is to be chosen suitably.
The purpose of this is that when we will vary $a$ over a small interval
and fix $b$ to be the midpoint of this interval,
we will obtain variables $\tilde{V}_n^+(a, b)$ that are defined with
the same drift,
\[
- \biggl(\frac{\vert f'(g(b)) \vert}{2 (L'(g(b)))^2} - 2\epsilon_n \biggr) t^2 ,
\]
and the Browian motion $W_{g(a)}$ only depending on $a$.
The case of $\tilde{V}_n^-(a, b)$ is similar.

For $0\leq u<v\leq1$, and $\alpha_n,\beta_n$ satisfying the conditions
of Theorem~\ref{theolawUn},
we obtain
\[
S_n\le \sup_{a\in[f(v-\beta_n),f(u+\alpha_n)]} \frac
{A(b(a))}{L'(g(b(a)))}\bigl\vert\tilde
V_n^+\bigl(a,b(a)\bigr)\bigr\vert\vee O_p(1) +
o_p(\log n)^{-2/3},
\]
and
\[
S_n\vee O_p(1)\ge\sup_{a\in[f(v-\beta_n),f(u+\alpha_n)]} \frac
{A(b(a))}{L'(g(b(a)))}
\bigl\vert\tilde V_n^-\bigl(a,b(a)\bigr)\bigr\vert+ o_p(\log
n)^{-2/3}
\]
for any $b(a)\in\R$ that satisfies $|a-b(a)|\le n^{-1/3}(\log n)^2$,
where $S_n$ is defined by~\eqref{defSn} and $\epsilon_n=1/\log n$
in~\eqref{defVntilde-} and~\eqref{defVntilde+}.
See Lemma~6.8 in~\cite{durotkulikovlopuhaa2012}.

Note that in order to obtain the above approximations, we use the
following lemma, which is a variation on Lemma~2.1 in~\cite
{kulikov-lopuhaa2006}. Although very simple, it turns out to be a very
useful tool to compare locations of maxima.
%
\begin{lem}\label{lemapproxargmax}
Let $I\subset\R$ be an interval. Let $g$ and $Z$ be real valued
functions defined on $I$ such that there exists $\gamma>0$ with
%
\[
g(u)<g(v)\qquad\mbox{for all } u,v \mbox{ such that }|u|>|v|+
\gamma.
\]
Assume that both $\sup_{u\in I} Z(u)$ and $\sup_{u\in I} \{Z(u)+g(u)\}$
are achieved.
Denoting by $\argmax$ an arbitrary point where the maximum is achieved,
we have
\[
\Bigl\llvert \mathop{\argmax}_{u\in I}\bigl\{Z(u)+g(u)\bigr\}\Bigr\rrvert \le\Bigl
\llvert \mathop{\argmax}_{u\in I}\bigl\{ Z(u)\bigr\}\Bigr\rrvert +\gamma.
\]
\end{lem}

\begin{pf}
Suppose the maximum of $Z$ is achieved at $v\in I$, so that $Z(u)\le
Z(v)$ for all $u\in I$.
It is assumed that for all $u\in I$ such that $|u|>|v|+\gamma,$ we have
$g(u)<g(v)$. Therefore,
\[
Z(u)+g(u)<Z(v)+g(v)
\]
for all $u\in I$ such that $|u|>|v|+\gamma$. It follows that the
maximum of $Z+g$ cannot be achieved at such a point $u$, which means that
\[
\Bigl\llvert \mathop{\argmax}_{u\in I}\bigl\{Z(u)+g(u)\bigr\}\Bigr\rrvert \le|v|+
\gamma.
\]
This completes the proof by definition of $v$.
\end{pf}

To relate the suprema of $\tilde{V}_n^+$ and $\tilde{V}_n^-$ with
maxima of independent random variables,
we will partition the interval $[f(v-\beta_n),f(u+\alpha_n)]$ into
a~union of disjoint intervals $A_i$ and $B_i$ of alternating length,
and a remainder interval~$R_n$, in such a way that the length of the
small blocks $A_i$ is
%
\begin{equation}
\label{deflengthln} l_n=\frac{2\Vert f' \Vert_\infty}{\inf_{t\in[0,1]}L'(t)}n^{-1/3} \log n,
\end{equation}
and the length of the big blocks $B_i$ is
$L_n = 2n^{-1/3} (\log n)^2$.
More precisely, for $i=1,2,\ldots,K_n$, where
%
\begin{equation}
\label{defKn} K_n= \biggl[\frac{f(u+\alpha_n)-f(v-\beta_n)}{l_n+L_n} \biggr]-1,
\end{equation}
let
%
\begin{eqnarray}
\label{defblocks} A_i &=& \bigl[f(v-\beta_n)+(i-1)
(l_n+L_n),  f(v-\beta_n)+il_n+(i-1)L_n
\bigr],
\nonumber
\\[-8pt]
\\[-8pt]
\nonumber
B_i &=&\bigl[f(v-\beta_n)+il_n+(i-1)L_n,
 f(v-\beta_n)+i(l_n+L_n)\bigr],
\end{eqnarray}
and let $R_n=[f(v-\beta_n)+K_n(l_n+L_n),f(u+\alpha_n)]$, so that
$l_n+L_n\leq|R_n|<2(l_n+L_n)$ and
%
\begin{equation}
\label{defdecomposition} \bigl[f(v-\beta_n),f(u+
\alpha_n)\bigr]= \Biggl(\bigcup_{i=1}^{K_n}
A_i \Biggr)\cup \Biggl(\bigcup_{i=1}^{K_n}
B_i \Biggr)\cup R_n.
\end{equation}
Now, suppose that $0\leq u<v\leq1$, and $\alpha_n,\beta_n$ satisfy the
conditions of Theorem~\ref{theolawUn}
and let $(\zeta_i)_{i\in\N}$ be a sequence of independent processes,
all distributed like $\zeta$ given in \eqref{defzeta}.
Then,
using scaling properties of the Brownian motion, we can build (possibly
dependent) copies~$(\zeta_j^{(1)})$,~$(\zeta_j^{(2)})$ of $(\zeta_i)_{i\in\N}$, such that
%
\begin{equation}\qquad
\label{eqUnxi+} S_B \leq \frac{S_n\vee O_p(1)}{1+O(1/\log n)}
\\
\leq S_B^{(1)} \vee S_A^{(2)} \vee
o_p(\log n)^{1/3} + o_p(\log
n)^{-2/3},
\end{equation}
where
\begin{eqnarray*}
S_B &\stackrel{d} {=}& \max_{1\le i\le K_n} \sup_{c\in[0,\Delta_{in}] }\bigl|
\zeta_i(c)\bigr| \quad\mbox{and}\quad S_B^{(1)} \stackrel{d}
{=} \max_{1\le i\le K_n} \sup_{c\in[0,\Delta_{in}] }\bigl| \zeta^{(1)}_i(c)\bigr|,
\\
S_A^{(2)} &\stackrel{d} {=} &\max_{2\le i\le K_n}
\sup_{c\in[0,\delta_{in}] }\bigl|\zeta^{(2)}_i(c)\bigr|,
\end{eqnarray*}
with $K_n$ defined in~\eqref{defKn} and where uniformly in $i$,
\[
\Delta_{in}=\bigl(1+o(1)\bigr) (\log n)^2\biggl\llvert
\frac{L'(g(b_i))f'(g(b_i))}{2}\biggr\rrvert^{-1/3},
\]
and $0\le\delta_{in}\le C\log n$, for some $C>0$, where $b_i$ denotes
the midpoint of the interval $B_i$ defined in~\eqref{defblocks}.
See Lemma~6.9 in~\cite{durotkulikovlopuhaa2012} for a
rigorous proof of~\eqref{eqUnxi+}.
The fact that $i\geq2$ in the definition of $S_A^{(2)}$ is due to the
fact that the first small block~$A_1$
has to be treated separately.

At this stage, we need a precise control of the tail probabilities of
the supremum of the limiting process $\zeta$ over increasing intervals.
Specifically, in Supplement A of~\cite{durotkulikovlopuhaa2012}, we obtain
the following slight variation on Theorem~1.1 in~\cite
{hooghiemstra-lopuhaa1998}.
Suppose $\delta_n\to\infty$, $\tau_n\to0$ and $u_n\to\infty$, in such
way that
$u_n/\delta_n\to0$,
$\delta_n\mu(u_n)/\tau_n\to1$,
and
$\log(\tau_n)/\delta_n^3\to0$.
Then
%
\begin{equation}
\label{eqextremallimitzeta} \biggl\llvert \frac{\log\prob (\sup_{c\in[0,\delta_n]}|\zeta(c)|\le u_n
)}{-2\tau_n}-1 \biggr\rrvert
\to0.
\end{equation}
See Lemma~6.3 in~\cite{durotkulikovlopuhaa2012}
for a rigorous proof.

We are then in the position to establish Theorem~\ref{theolawUn} and
Corollary~\ref{corlawUn}.

\begin{pf*}{Proof of Theorem~\protect\ref{theolawUn}}
Let $(u_n)_n$ be a sequence such that $u_n\to\infty$ in such a way that
%
\begin{equation}
\label{eqsequenceun} n^{1/3}\mu(u_n)\to\tau>0,
\end{equation}
where $\mu$ is the density of $\zeta(0)$.
We will bound
$\P(S_n\leq u_n)$, where $S_n$ is defined by~\eqref{defSn},
from above and below by means of~\eqref{eqUnxi+}.
Write
\begin{eqnarray*}
S_1 &= &\max_{1\le i\le K_n} \sup_{c\in[0,\Delta_{in}] }\bigl|
\zeta^{(1)}_i(c)\bigr|,
\\
S_2 &=& \max_{2\le i\le K_n}\sup_{c\in[0,\delta_{in}] }\bigl|\zeta^{(2)}_i(c)\bigr|.
\end{eqnarray*}
Then, according to \eqref{eqUnxi+}
\[
\P(S_n\leq u_n) \ge \prob \bigl( \bigl(1+O(1/\log n)
\bigr) \{ S_1\vee S_2\vee Q_n \}+
R_n\le u_n \bigr),
\]
where $Q_n=o_p(\log n)^{1/3}$ and $R_n=o_p(\log n)^{-2/3}$.
Define the event $E_n=\{(\log n)^{2/3}|R_n|\le1\}$, then $\prob
(E_n^c)\to0$,
so that 
\begin{eqnarray*}
\P(S_n\leq u_n) &\ge& \prob ( S_1\vee
S_2\vee Q_n\le v_n ) + o(1)
\\
&=& \prob ( S_1\le v_n,S_2\le
v_n,Q_n\le v_n ) +o(1),
\end{eqnarray*}
where
\[
v_n = \frac{u_n-(\log n)^{-2/3}}{1+O((\log n)^{-1})} \sim u_n-(\log
n)^{-2/3} \qquad\mbox{as $n\to\infty$.}
\]
From \eqref{eqexpansionmu} and \eqref{eqsequenceun}, it is easily
verified that
$u_n$ is of order $(\log n)^{1/3}$ [see also the expansion~(\ref{eqexpansionun}) below] and that
%
\begin{equation}
\label{eqsequencevn} n^{1/3}\mu(v_n)\to\tau.
\end{equation}
Therefore, since $\prob(Q_n\le v_n)\to1$,
we have
\[
\P(S_n\leq u_n)\ge\prob (S_1\le
v_n,S_2\le v_n )+o(1).
\]
We will investigate $\prob (S_1\le v_n )$ and
$\prob (S_2\le v_n )$ separately.

Since the processes $\zeta_i^{(1)}$ are independent copies of $\zeta$,
\[
\prob (S_1\le v_n ) = \prod
_{i=1}^{K_n} \prob \Bigl( \sup_{c\in[0,\Delta_{in}] }\bigl|
\zeta(c)\bigr|\le v_n \Bigr).
\]
For each $i=1,2,\ldots,K_n$ fixed, we apply~\eqref{eqextremallimitzeta},
with
\[
\Delta_{in}=\bigl(1+o(1)\bigr) (\log n)^2\biggl\llvert
\frac{L'(g(b_i))f'(g(b_i))}{2}\biggr\rrvert^{-1/3},
\]
which is of the order $(\log n)^2$ uniformly in $i$, and $\tau_{in}=\tau\Delta_{in} n^{-1/3}$,
where the $b_i$ are the midpoints of the $K_n$ big blocks $B_i$.
The $b_i$ are equidistant at distance
$l_n+L_n=2n^{-1/3}(\log n)^2(1+O(\log n)^{-1})$.
Since $\tau_{in}\to0$ uniformly\vadjust{\goodbreak} in $i$ and $v_n$ is of order $(\log
n)^{1/3}$, we conclude that
\begin{eqnarray*}
&&\prod_{i=1}^{K_n} \prob \Bigl(
\sup_{c\in[0,\Delta_{in}] }\bigl| \zeta(c)\bigr|\le v_n \Bigr)\\
&&\qquad = \prod
_{i=1}^{K_n} \exp \bigl( -2\tau_{in}
\bigl(1+o(1)\bigr) \bigr),
\end{eqnarray*}
where the small $o$-term is uniform in $i$.
Therefore,
\begin{eqnarray*}
& &\prod_{i=1}^{K_n} \prob \Bigl(
\sup_{c\in[0,\Delta_{in}] }\bigl| \zeta(c)\bigr|\le v_n \Bigr)
\\
&&\qquad= \exp \Biggl\{ -2\bigl(1+o(1)\bigr)\tau \sum_{i=1}^{K_n}
\frac{2n^{-1/3}(\log n)^2}{\llvert 4L'(g(b_i))f'(g(b_i))\rrvert^{1/3}} \Biggr\}
\\
&&\qquad= \exp \biggl\{ -2\tau\int_{f(v)}^{f(u)}
\frac{1}{\llvert 4L'(g(b))f'(g(b))\rrvert^{1/3}} \,\mathrm{d}b \biggr\}+o(1)
\\
&&\qquad= \exp \biggl\{ -2\tau\int_u^v
\frac{|f'(t)|^{2/3}}{(4L'(t))^{1/3}} \,\mathrm{d}t \biggr\}+o(1).
\end{eqnarray*}
It follows that
\[
\prob (S_1\le v_n ) \to \exp \biggl\{ -2\tau\int
_u^v \frac{|f'(t)|^{2/3}}{(4L'(t))^{1/3}} \,\mathrm{d}t \biggr\}.
\]
The probability $\prob (S_2\le v_n )$ can be treated in the
same way:
\begin{eqnarray*}
\prob (S_2\le v_n ) &=& \prod
_{i=1}^{K_n} \prob \Bigl( \sup_{c\in[0,\delta_{in}] }\bigl|
\zeta(c)\bigr|\le v_n \Bigr)
\\
&=& \exp \Biggl\{ -2\bigl(1+o(1)\bigr)\tau \sum_{i=1}^{K_n}
\delta_{in}n^{-1/3} \Biggr\} \to1,
\end{eqnarray*}
since, according to~\eqref{eqUnxi+}
and \eqref{defKn},
\[
\sum_{i=1}^{K_n}\delta_{in}n^{1/3}
\le Cn^{-1/3}K_n\log n=O(\log n)^{-1}.
\]
This yields that
\[
\liminf_{n\to\infty} \prob (S_n\le u_n ) \ge \exp
\biggl\{ -2\tau\int_u^v \frac{|f'(t)|^{2/3}}{(4L'(t))^{1/3}}
\,\mathrm{d}t \biggr\}.
\]
Similarly, with~\eqref{eqUnxi+},
\[
\P(S_n\leq u_n)\leq\prob \Bigl( \max_{1\le i\le K_n}
\sup_{c\in[0,\Delta_{in}] }\bigl| \zeta_i(c)\bigr|\le v_n \Bigr) + o(1),\vadjust{\goodbreak}
\]
where $v_n$ satisfies \eqref{eqsequencevn}.
This probability can be treated completely similar to $\prob(S_1\le
v_n)$, so that
\[
\limsup_{n\to\infty} \prob (S_n\le u_n ) \le \exp
\biggl\{ -2\tau\int_u^v \frac{|f'(t)|^{2/3}}{(4L'(t))^{1/3}}
\,\mathrm{d}t \biggr\}.
\]
This proves the theorem.
\end{pf*}

\begin{pf*}{Proof of Corollary~\protect\ref{corlawUn}}
Let $(u_n)_n$ be a sequence such that $u_n\to\infty$ in such a way that
$n^{1/3}\mu(u_n)\to\tau>0$, as $n\to\infty$.
Taking logarithms in~\eqref{eqexpansionmu}, we conclude that
$(u_n)_n$ should satisfy
%
\begin{equation}
\label{eqexpansionlogmu} \frac13\log n+\log u_n-
\frac23u_n^3-\kappa u_n = \log
\frac{\tau}{2\lambda}+o(1)\qquad \mbox{as } n\to\infty.
\end{equation}
This means that $-2u_n^3/3$ is the dominating term, which should
compensate $(\log n)/3$.
Therefore, if we write $u_n=2^{-1/3}(\log n)^{1/3}+\delta_n$, where
$\delta_n=  o(\log n)^{1/3}$, and insert this in~\eqref{eqexpansionlogmu}, we obtain
%
\begin{eqnarray*}
&&\frac13\log n + \log \biggl\{ \biggl(\frac{\log n}{2} \biggr)^{1/3}+
\delta_n \biggr\}
\\
&&\quad{}- \frac23 \biggl\{\frac{\log n}{2} + 3 \biggl(\frac{\log n}{2}
\biggr)^{2/3}\delta_n + 3 \biggl(\frac{\log n}{2}
\biggr)^{1/3}\delta_n^2+ \delta_n^3
\biggr\}-\kappa \biggl(\frac{\log n}{2} \biggr)^{1/3}-\kappa
\delta_n\\
&&\qquad = \log\frac{\tau}{2\lambda}+o(1).
\end{eqnarray*}

Tedious, but straightforward computations first yield that $\delta_n\to
0$ and then that
\begin{eqnarray*}
\delta_n &=& -\frac{\kappa}{4^{1/3}}(\log n)^{-1/3} +
\frac{4^{1/3}}{6}(\log n)^{-2/3}\log\log n
\\
&&{} -(\log n)^{-2/3} \biggl[ \frac{\log\tau}{2^{1/3}}-\frac{\log(2\lambda)}{2^{1/3}} +
\frac{4^{1/3}}{6}\log2 \biggr] +o(\log n)^{-2/3}.
\end{eqnarray*}
%
If we put $\tau4^{-1/3}C_{f,L}=\mathrm{e}^{-x}$, or $-\log\tau=x+\log
C_{f,L}-(2\log2)/3$, this implies that
%
\begin{eqnarray}
\label{eqexpansionun}
u_n &= &\frac{1}{2^{1/3}}(\log
n)^{1/3} -\frac{\kappa}{4^{1/3}}(\log n)^{-1/3} +
\frac{4^{1/3}}{6}(\log n)^{-2/3}\log\log n
\nonumber
\\[-8pt]
\\[-8pt]
\nonumber
&&{} +(\log n)^{-2/3} \biggl( \frac{x+\log C_{f,L}}{2^{1/3}}+\frac{\log\lambda}{2^{1/3}} \biggr)
+ o(\log n)^{-2/3}.
\end{eqnarray}
If we also write $u_n=x/a_n+b_n+o(\log n)^{-2/3}$, with
\begin{eqnarray*}
a_n &=& 2^{1/3}(\log n)^{2/3},
\\
b_n &=& \frac{(\log n)^{1/3}}{2^{1/3}} - \frac{\kappa}{(4\log n)^{1/3}} + \frac{4^{1/3}\log\log n}{6(\log n)^{2/3}}
+ \frac{\log(\lambda C_{f,L})}{2^{1/3}(\log n)^{2/3}},
\end{eqnarray*}
then~\eqref{eqLimitsupU} is equivalent to
$\P \{
a_n(S_n-b_n)+o(1)\leq x \}\to\exp \{-\mathrm{e}^{-x} \}$.
Finally, it is easy to see that
\begin{eqnarray*}
a_n(S_n-b_n) &=& \log n \biggl\{ \biggl(
\frac{2}{\log n} \biggr)^{1/3}S_n- \biggl(
\frac{2}{\log n} \biggr)^{1/3} b_n \biggr\}
\\
&=& \log n \biggl\{ \biggl( \frac{2}{\log n} \biggr)^{1/3}S_n-
\mu_n \biggr\}.
\end{eqnarray*}
This proves Corollary~\ref{corlawUn}.
\end{pf*}

\section{\texorpdfstring{Proof of Theorems~\protect\ref{theoratefn} and \protect\ref{theolawfn}}
{Proof of Theorems 2.1 and 2.2}}
\label{secproofratelawfn}
We suppose in the sequel that assumptions (A1), (A2) and (A3) are fulfilled.
As before, $C$, $C_1$, $C_2,\ldots$ denote positive real numbers that
depend only on $q$, $C_q$, $f$, $L$, $C_0$,
and possibly also on $\holder$ under the additional assumption (A4).
It follows from the definition of~$\hat f_n$ that it can be
discontinuous only at the jump points of~$F_n$.
In particular, the number of jump points of $\hat f_n$ is finite.
In the sequel, we will denote this number by $N_n-1$ (note that $N_n\ge1$).
Moreover, we set $\tau_0=0$, $\tau_{N_n}=1$, and in the case where
$\hat f_n$ has at least one jump point, that is,
$N_n\ge2$, we denote by $\tau_1 < \cdots< \tau_{N_n-1}$ the ordered
jump points of $\hat f_n$.

To prove Theorems~\ref{theoratefn} and \ref{theolawfn},
we need a precise uniform bound on the spacings between consecutive
jump points of $\hat f_n$.
This is given by the following lemma.
%
\begin{lem}
\label{lemjump}
Assume \textup{(A1)} and \textup{(A2)}.
Then
%
\begin{equation}
\label{eqlemjump} \max_{1\le i\le N_n}|\tau_i-\tau_{i-1}|=O_p
\biggl(\frac{\log n}{n} \biggr)^{1/3}.
\end{equation}

\end{lem}
\begin{pf}
It follows from the definition of $\hat f_n$ and $\hat U_n$ that these
functions are nonincreasing left-continuous step functions with
finitely many jump points, and that the maximal length of the flat
parts of $\hat f_n$ is precisely the maximal height of the jumps of
$\hat U_n$. Therefore,
\[
\max_{1\le i\le N_n}|\tau_i-\tau_{i-1}| =
\sup_{a\in\R} \Bigl\vert \lim_{b\downarrow a}\hat U_n(b)-\hat
U_n(a) \Bigr\vert.
\]
Using the triangle inequality, it follows that
\[
\max_{1\le i\le N_n}|\tau_i-\tau_{i-1}| \le
\sup_{a\in\R} \Bigl\{ \Bigl\vert\lim_{b\downarrow a}\hat U_n(b)-g(a)
\Bigr\vert + \bigl\vert\hat U_n(a)-g(a) \bigr\vert \Bigr\}.
\]
But $g$ is continuous on $\R$, so that Theorem \ref{theorateUn}
implies that
\[
\max_{1\le i\le N_n}|\tau_i-\tau_{i-1}| \le 2
\sup_{a\in\R} \bigl\vert\hat U_n(a)-g(a) \bigr\vert= O_p
\biggl(\frac
{\log n}{n} \biggr)^{1/3},
\]
which completes the proof.
\end{pf}

%
\begin{remark}
\label{onthejumppts}
Lemma \ref{lemjump} together with the identity
$1=\sum_{i=1}^{N_n}(\tau_i-\tau_{i-1})$,
implies that
$1/N_n=O_p(n^{-1/3}(\log n)^{1/3})$.
This gives some idea about the order of magnitude of the number of
jumps of $\hat f_n$.
Further investigation is needed to obtain a sharp upper bound, and we
conjecture that it is of order~$n^{1/3}$.
This rate is also claimed in Theorem 3.1
in~\cite{groeneboom2011}.\vspace*{-2pt}
\end{remark}
We will also need a bound on the mean absolute error
between $\hat f_n$ and $f$. In Supplement C in~\cite
{durotkulikovlopuhaa2012}, we reprove Theorem 1 in~\cite{durot2007}
under slightly weaker assumptions; that is,
there exists $C>0$ such that
%
\begin{equation}
\label{eqEfn} \E\bigl|\hat f_n(t)-f(t)\bigr|\le Cn^{-1/3}
\end{equation}
for all $t\in[n^{-1/3},1-n^{-1/3}]$ and
%
\begin{equation}
\label{eqEfnbord} \E\bigl|\hat f_n(t)-f(t)\bigr|\le C \bigl[n\bigl(t\wedge(1-t)
\bigr) \bigr]^{-1/2}
\end{equation}
for all $t\in(0,n^{-1/3}]\cup[1-n^{-1/3},1)$.
See Lemma~6.10 in~\cite{durotkulikovlopuhaa2012}.

Note that the number of jump points of $\hat U_n$ is precisely the
number of flat parts of $\hat f_n$, that is $N_n$, and
denoting by $\gamma_1>\cdots>\gamma_{N_n}$ the jump points of~$\hat
U_n$, we have
%
\begin{equation}
\label{eqgammatau} \gamma_i=\hat f_n(
\tau_i)\quad\mbox{and}\quad \tau_i=\hat U_n(
\gamma_i)\qquad\mbox{for  } i=1,2,\ldots,N_n.
\end{equation}
We will show that in order to study the supremum of $|\widehat{f}_n-f|$
over an interval,
we can restrict ourselves to the situation where the boundaries of the
interval are jump points of $\widehat{f}_n$
and where the values of $\hat f_n$ stay in $(f(1), f(0))$.
Indeed, in order to relate the supremum of $|\widehat{f}_n-f|$ to the
supremum of $|U_n-g|$,
we need to employ the identity $\gamma_i=f(g(\gamma_i))$, for
$\gamma_i=\widehat{f}_n(\tau_i)$, so we need to make sure that $\widehat
{f}_n(\tau_i)\in(f(1),f(0))$.
To this end, define
for any $t\in(0,1)$
%
\begin{equation}
\label{defi1} i_1(t)=\min\bigl\{i\in\{1,2,\ldots,N_n\}
\mbox{ such that } \tau_i\ge t\bigr\}
\end{equation}
and
%
\begin{equation}
\label{defi2} i_2(t)=\max\bigl\{i\in\{0,1,\ldots,N_n-1\}
\mbox{ such that } \tau_i<1-t\bigr\}.
\end{equation}
For any $t$ such that $n^{1/3}t\to\infty$ and $n^{1/3}(1-t)\to\infty$,
we establish the order of the difference with neighboring points of
jump of $\widehat{f}_n$, that is,
%
\begin{equation}
\label{eqlemtau1} \tau_i=t+O_p\bigl(n^{-1/3}
\bigr)
\end{equation}
for $i=i_1(t)-1,i_1(t),i_1(t)+1$, and similarly for $1-t$,
%
\begin{equation}
\label{eqlemtau2} \tau_i=1-t+O_p\bigl(n^{-1/3}
\bigr)
\end{equation}
for $i=i_2(t)-1,i_2(t),i_2(t)+1$.
See Lemma~6.11 in~\cite{durotkulikovlopuhaa2012}.
Note that if there are no jumps on the interval $[s,1-t)$, then $\tau_{i_1(s)}>\tau_{i_2(t)}$.
This may happen if the length $1-t-s$ of the interval tends to zero too fast.
However, if
%
\begin{equation}
\label{eqcondEn} n^{1/3}s\to\infty, \qquad n^{1/3}t\to\infty
\end{equation}
and
%
\begin{equation}
\label{eqcond73} n^{1/3}(1-t-s)\to\infty,\vadjust{\goodbreak}
\end{equation}
then
%
\begin{equation}
\label{eqlemorderingi1i2} \P(s\leq\tau_{i_1(s)}\leq\tau_{i_2(t)}<1-t)
\to1.
\end{equation}
See Lemma~6.12 in~\cite{durotkulikovlopuhaa2012}.
According to Lemma~6.13 in~\cite{durotkulikovlopuhaa2012},
%
\begin{eqnarray}
\label{eqlemEn}
\P \bigl( \gamma_i<f(0) \mbox{ for all }i\geq
i_1(s) \bigr) &\to&1,
\nonumber
\\[-8pt]
\\[-8pt]
\nonumber
\P \bigl( \gamma_i>f(1) \mbox{ for all }i\leq i_2(t)
\bigr) &\to&1,
\end{eqnarray}
whenever~\eqref{eqcond73} holds,
which ensures that $\widehat{f}_n(\tau_i)\in(f(1),f(0))$ simultaneously
for various $i$'s, with probability tending to one.

We are then in the position to prove Theorem~\ref{theoratefn}.

\begin{pf*}{Proof of Theorem~\protect\ref{theoratefn}}
First, we establish the result for sequences $\alpha_{n}=s_n$ and $\beta_{n}=t_n$
that satisfy~\eqref{eqcondEn} and~\eqref{eqcond73}.
For the sake of brevity, write $i_1=i_1(s_n)$ and $i_2=i_2(t_n)$.
Define the event
%
\begin{eqnarray}
\label{eqE}
E_n &= &\{ s_n\le\tau_{i_1}\le
\tau_{i_2}<1-t_n \}
\nonumber
\\[-8pt]
\\[-8pt]
\nonumber
&& {}\cap \bigl\{\gamma_i\in\bigl(f(1),f(0)\bigr) \mbox{ for all
}i=i_1,\ldots, i_2 \bigr\}.
\end{eqnarray}
Then according to~\eqref{eqlemorderingi1i2} and~\eqref{eqlemEn},
we have $\P(E_n)\to1$,
so we can
restrict ourselves to the event $E_n$.
We have
\begin{eqnarray*}
&&\sup_{u\in(s_n,1-t_n]}\bigl|\hat f_n(u)-f(u)\bigr|
\\
&&\qquad \le \max_{i=i_1,\ldots,i_2}\sup_{u\in(\tau_{i-1},\tau_i]}\bigl|\hat f_n(u)-f(u)\bigr| +
\sup_{u\in(\tau_{i_2},1-t_n]}\bigl|\hat f_n(u)-f(u)\bigr|.
\end{eqnarray*}
Recall that $\hat f_n$ is constant on every interval $(\tau_{i-1},\tau_i]$, for $i=1,2,\ldots,N_n-1$.
Moreover, $f'$ is bounded. Using the triangle inequality, it follows that
\begin{eqnarray*}
\sup_{u\in(\tau_{i-1},\tau_i]}\bigl|\hat f_n(u)-f(u)\bigr| &=& \sup_{u\in(\tau_{i-1},\tau_i]}\bigl|\hat
f_n(\tau_i)-f(u)\bigr|
\\
&\le &\bigl|\hat f_n(\tau_i)-f(\tau_i)\bigr| +
\bigl\|f'\bigr\|_\infty|\tau_{i-1}-\tau_i|
\end{eqnarray*}
for all $i=1,2,\ldots,N_n-1$ and
\[
\sup_{u\in(\tau_{i_2},1-t_n]}\bigl|\hat f_n(u)-f(u)\bigr| \le\bigl |\hat f_n(1-t_n)-f(1-t_n)\bigr|+
\bigl\|f'\bigr\|_\infty|\tau_{i_2}-\tau_{i_2+1}|.
\]
From~\eqref{eqcondEn} and~\eqref{eqcond73},
we have $1-t_n\in[n^{-1/3},1-n^{-1/3}]$, for large enough $n$, so~(\ref
{eqEfn}) ensures that
$\hat f_n(1-t_n)-f(1-t_n)=O_p(n^{-1/3})$.
Using~\eqref{eqlemjump} and~(\ref{eqgammatau}), it follows that
\begin{eqnarray*}
\sup_{u\in(s_n,1-t_n]}\bigl|\hat f_n(u)-f(u)\bigr| &\le& \max_{i=i_1,\ldots,i_2}\bigl|\hat
f_n(\tau_i)-f(\tau_i)\bigr| +O_p
\biggl(\frac{\log
n}{n} \biggr)^{1/3}
\\[-2pt]
&=& \max_{i=i_1,\ldots,i_2}\bigl|\gamma_i-f\circ\hat U_n(
\gamma_i)\bigr| +O_p \biggl(\frac{\log n}{n}
\biggr)^{1/3}.
\end{eqnarray*}
On the event $E_n$, we have $\gamma_i=f\circ g(\gamma_i)$, for all
$i=i_1,\ldots,i_2$, and therefore
\begin{eqnarray*}
\sup_{u\in(s_n,1-t_n]}\bigl|\hat f_n(u)-f(u)\bigr| &\le &\bigl\|f'
\bigr\|_\infty\max_{i=i_1,\ldots,i_2}\bigl|g(\gamma_i)-\hat
U_n(\gamma_i)\bigr| +O_p \biggl(
\frac{\log n}{n} \biggr)^{1/3}
\\[-2pt]
&\le &\bigl\|f'\bigr\|_\infty\sup_{a\in\R}\bigl|\hat
U_n(a)-g(a)\bigr|+O_p \biggl(\frac{\log
n}{n}
\biggr)^{1/3}.
\end{eqnarray*}
Theorem~\ref{theoratefn}, with $\alpha_n=s_n$ and $\beta_n=t_n$
satisfying~\eqref{eqcondEn} and~\eqref{eqcond73}
now follows from Theorem~\ref{theorateUn}.

It remains to extend the result to more general sequences $\alpha_n$
and $\beta_n$.
For this purpose, define $s_n=n^{-1/3}(\log n)^{1/6}$.
In view of the foregoing results, we know that
%
\begin{equation}
\label{eqsupan} \sup_{t\in(s_n,1-s_n]}\bigl|\hat f_n(t)-f(t)\bigr|=O_p
\biggl(\frac{\log n}{n} \biggr)^{1/3}.
\end{equation}
Suppose $\alpha_n$ and $\beta_n$ satisfy (\ref{eqalphanratefn}).
Let us notice that
$\sup_{t\in(\alpha_n,1-\beta_n]}|\hat f_n(t)-f(t)|$
decreases when either $\alpha_n$ or $\beta_n$ increases, so that we can
restrict our attention to small values of $\alpha_n$ and $\beta_n$.
Without loss of generality we may assume that $\alpha_n\leq
n^{-1/3}\leq s_n$ and $\beta_n\leq n^{-1/3}$.

We then use the following property of nonincreasing functions $h_1$ and $h_2$
on an interval $[a,b]$:
%
\begin{eqnarray}
\label{eqlemmonotone} &&\sup_{t\in[a,b]}\bigl|h_1(t)-h_2(t)\bigr|
\nonumber
\\[-9pt]
\\[-9pt]
\nonumber
&&\qquad\le\bigl|h_1(a)-h_2(a)\bigr|\vee \bigl|h_1(b)-h_2(b)\bigr|+\bigl|h_2(a)-h_2(b)\bigr|.
\end{eqnarray}
See Lemma~6.1 in~\cite{durotkulikovlopuhaa2012}.
Since $\hat f_n$ and $f$ are both nonincreasing,
according to~\eqref{eqlemmonotone}, we have
\begin{eqnarray*}
&& \sup_{t\in(\alpha_n,s_n]}\bigl|\hat f_n(t)-f(t)\bigr|
\\
&&\qquad \le \bigl|\hat f_n(\alpha_n)-f(\alpha_n) \bigr|\vee
\bigl|f(s_n)-\hat f_n(s_n) \bigr|+\bigl\|f'
\bigr\|_\infty(s_n-\alpha_n).
\end{eqnarray*}
Because $s_n\in[n^{-1/3}, 1-n^{-1/3}]$,
it follows from~\eqref{eqEfn} and~\eqref{eqEfnbord} that
$f(s_n)-\hat f_n(s_n)=O_p(n^{-1/3})$
and
$\hat f_n(\alpha_n)-f(\alpha_n)=O_p((n\alpha_n)^{-1/2})$,
which is of the order $O_p(n^{-1/3}(\log n)^{1/3})$,
as we have assumed that $\alpha_n\ge K_1n^{-1/3}(\log n)^{-2/3}$.
We conclude
\[
\sup_{t\in(\alpha_n,s_n]}\bigl|\hat f_n(t)-f(t)\bigr|=O_p \biggl(
\frac{\log
n}{n} \biggr)^{1/3}.
\]
Similarly, we obtain
\[
\sup_{t\in(1-s_n,1-\beta_n]}\bigl|\hat f_n(t)-f(t)\bigr|=O_p \biggl(
\frac{\log
n}{n} \biggr)^{1/3
}
\]
and therefore,
\[
\sup_{t\in(\alpha_n,1-\beta_n]}\bigl|\hat f_n(t)-f(t)\bigr|=\sup_{t\in
(s_n,1-s_n]}\bigl|\hat
f_n(t)-f(t)\bigr|\vee O_p \biggl(\frac{\log n}{n}
\biggr)^{1/3}.
\]
Theorem \ref{theoratefn} now follows from (\ref{eqsupan}).
\end{pf*}

To prove Theorem~\ref{theolawfn}, similarly to the proof of
Theorem~\ref{theoratefn},
we first establish the result for sequences $s_n=u+\alpha_n$ and
$t_n=v-\beta_n$
satisfying~\eqref{eqcondEn} and~\eqref{eqcond73},
and then extend the result to more general sequences.
The first step is to prove that the behavior of supremum over the
interval $(s_n,1-t_n]$
is dominated by that of the largest interval between
two jump points of~$\widehat{f}_n$ contained in $(s_n,1-t_n]$.
For this task, we make use of the notation $\tau_i$, $\gamma_i$, $i_1$
and $i_2$ as introduced in~\eqref{eqgammatau},~\eqref{defi1}
and~\eqref{defi2}, and for $t\in[0,1]$, we define the normalizing function
%
\begin{equation}
\label{defB} B(t)=\bigl(4\bigl|f'(t)\bigr|L'(t)
\bigr)^{-1/3}.
\end{equation}
It is easy to see that under assumptions (A1), (A2) and (A4),
there exists $C_0>0$ and $\holder\in(0,1]$ such that
%
\begin{equation}
\label{eqABholder} \bigl|A(u)-A(v)\bigr|\le C_0|u-v|^\holder\quad \mbox{and}\quad
\bigl|B(u)-B(v)\bigr|\le C_0|u-v|^\holder
\end{equation}
for all $u,v\in[0,1]$, where $A$ is given by \eqref{defA}.
Recall that by convention, the supremum of an empty set is equal to zero.

For $s,t$ that satisfy conditions~\eqref{eqcond73} and~\eqref
{eqcondEn},
we first obtain
%
\begin{eqnarray}
\label{eqlemi1-i2}
&&\sup_{u\in(s,1-t]} B(u)\bigl\llvert \widehat{f}_n(u)-f(u)
\bigr\rrvert
\nonumber
\\[-8pt]
\\[-8pt]
\nonumber
&&\qquad = \sup_{u\in(\tau_{i_1(s)},\tau_{i_2(t)}]} B(u)\bigl\llvert \widehat{f}_n(u)-f(u)\bigr
\rrvert \vee O_p \bigl(n^{-1/3} \bigr).
\end{eqnarray}
See Lemma~6.14 in Supplement C in~\cite{durotkulikovlopuhaa2012}.
We are then able to make the connection between $\hat U_n$ and $\hat f_n$.
%
\begin{lem}
\label{lemslope-inverse}
Assume \textup{(A1), (A2), (A3)} and \textup{(A4)}.
Let $0<s<1-t<1$, possibly depending on $n$, such that $s,t$ satisfy
conditions~\eqref{eqcondEn}
and~\eqref{eqcond73}.
Then
\begin{eqnarray*}
&&\sup_{u\in(s,1-t]} B(u)\bigl\llvert \widehat{f}_n(u)-f(u)\bigr
\rrvert
\\
& &\qquad= \sup_{a\in[f(1-t),f(s)]} A(a)\bigl|\widehat{U}_n(a)-g(a)\bigr| +
O_p \biggl(\frac{\log n}{n} \biggr)^{({\holder+1})/{3}}
\end{eqnarray*}
for some $\holder\in(0,1]$.
\end{lem}

\begin{pf}
Again write $i_1=i_1(s)$ and $i_2=i_2(t)$.
We first decompose the supremum into maxima of suprema taken
over\vadjust{\goodbreak}
intervals between
succeeding jump points of $\widehat{f}_n$:
\[
\sup_{u\in(\tau_{i_1},\tau_{i_2}]} B(u)\bigl\llvert \widehat{f}_n(u)-f(u)\bigr
\rrvert = \max_{i_1+1\le i\le i_2} \sup_{u\in(\tau_{i-1},\tau_{i}]} B(u)\bigl\llvert
\widehat{f}_n(u)-f(u)\bigr\rrvert .
\]
Then, by Theorem~\ref{theoratefn} and~\eqref{eqlemjump},
we have that
%
\[
\sup_{u\in(\tau_{i_1},\tau_{i_2}]} \bigl\llvert \widehat{f}_n(u)-f(u)
\bigr\rrvert \le \sup_{u\in(s,1-t]} \bigl\llvert \widehat{f}_n(u)-f(u)
\bigr\rrvert =O_p \biggl(\frac{\log n}{n} \biggr)^{1/3}.
\]
Thus, we obtain by means of~(\ref{eqABholder}) and the triangle
inequality that
\begin{eqnarray*}
&&\sup_{u\in(\tau_{i_1},\tau_{i_2}]} B(u)\bigl\llvert \widehat{f}_n(u)-f(u)\bigr
\rrvert
\\
&&\qquad = \max_{i_1+1\le i\le i_2} B(\tau_i) \sup_{u\in(\tau_{i-1},\tau_{i}]} \bigl
\llvert \widehat{f}_n(u)-f(u)\bigr\rrvert + O_p \biggl(
\frac{\log n}{n} \biggr)^{{(\holder+1)}/{3}}.
\end{eqnarray*}
By monotonicity of $f$, {we have} for all $i_1+1\le i\le i_2$,
\[
\sup_{u\in(\tau_{i-1},\tau_{i}]} \bigl\llvert \widehat{f}_n(
\tau_i)-f(u)\bigr\rrvert = \bigl\llvert \widehat{f}_n(
\tau_i)-f(\tau_{i})\bigr\rrvert \vee \bigl\llvert
\widehat{f}_n(\tau_i)-f(\tau_{i-1})\bigr\rrvert
.
\]
{Hence, with (\ref{eqgammatau})} we arrive at
\begin{eqnarray*}
&& {\sup_{u\in(\tau_{i_1},\tau_{i_2}]}} B(u)\bigl\llvert \widehat{f}_n(u)-f(u)\bigr
\rrvert
\\
&&\qquad= \max_{i_1+1\le i\le i_2}B(\tau_i) \bigl\{ \bigl\llvert
\gamma_i-f(\tau_{i})\bigr\rrvert \vee\bigl\llvert
\gamma_i-f(\tau_{i-1})\bigr\rrvert \bigr\} {+
O_p \biggl(\frac{\log n}{n} \biggr)^{{(\holder+1)}/{3}}}.
\end{eqnarray*}
On the event $E_n$ of~\eqref{eqE}, we can write $\gamma_i=f(g(\gamma_i))$ for all $i=i_1+1,\ldots,i_2,$
which, in view of (\ref{eqgammatau}), implies that
\begin{eqnarray*}
\bigl|\gamma_i-f(\tau_{i})\bigr| &=& \bigl|g(\gamma_i)-
\widehat{U}_n(\gamma_i)\bigr|\cdot\bigl|f'(
\theta_{i1})\bigr|,
\\
\bigl|\gamma_i-f(\tau_{i-1})\bigr| &=& \bigl|g(\gamma_i)-
\widehat{U}_n(\gamma_{i-1})\bigr|\cdot\bigl|f'(
\theta_{i2})\bigr|
\end{eqnarray*}
for some $\theta_{i1}$ between $g(\gamma_i)$ and $\widehat{U}_n(\gamma_i)$,
{and $\theta_{i2}$} between $g(\gamma_i)$ and $\widehat{U}_n(\gamma_{i-1})$.
By~\eqref{eqlemEn}, Theorem~\ref{theorateUn} and~\eqref{eqA4},
it follows that
%
\begin{eqnarray}
\label{eqgi-fti}\bigl |\gamma_i-f(\tau_{i})\bigr| 
{=} \bigl|g(\gamma_i)-
\widehat{U}_n(\gamma_i)\bigr|\cdot\bigl|f'\bigl(g(
\gamma_i)\bigr)\bigr|+O_p \biggl(\frac{\log n}{n}
\biggr)^{{(\holder+1)}/{3}}.
\end{eqnarray}
By \eqref{eqgammatau}, \eqref{eqlemjump} and Theorem \ref{theorateUn}, we have that
%
\begin{eqnarray}
\label{eqgi-Uni-1} \max_{i_1+1\le i\le i_2}\bigl|g(\gamma_i)-
\widehat{U}_n(\gamma_{i-1})\bigr| &=& \max_{i_1+1\le i\le i_2}\bigl|g(
\gamma_i)-\widehat{U}_n(\gamma_{i})+
\tau_{i}-\tau_{i-1}\bigr|
\nonumber\\
&\le&\sup_{a\in\R}\bigl|g(a)-\widehat{U}_n(a)\bigr|+O_p
\biggl(\frac{\log
n}{n} \biggr)^{1/3}
\\
&=& O_p \biggl(\frac{\log n}{n} \biggr)^{1/3},\nonumber
\end{eqnarray}
so that similarly as above,
\[
\bigl|\gamma_i-f(\tau_{i-1})\bigr| = \bigl|g(\gamma_i)-
\widehat{U}_n(\gamma_{i-1})\bigr|\cdot\bigl|f'\bigl(g(
\gamma_i)\bigr)\bigr| + O_p \biggl(\frac{\log n}{n}
\biggr)^{{(\holder+1)}/{3}}.
\]
It follows {that}
\begin{eqnarray*}
& &{\sup_{u\in(\tau_{i_1},\tau_{i_2}]}} B(u)\bigl\llvert \widehat{f}_n(u)-f(u)\bigr
\rrvert
\\[-2pt]
&&\qquad = \max_{i_1+1\le i\le i_2} B(\tau_i)\bigl|f'\bigl(g(
\gamma_i)\bigr)\bigr| \bigl\{ \bigl|g(\gamma_i)-
\widehat{U}_n(\gamma_i)\bigr| \vee \bigl|g(\gamma_i)-
\widehat{U}_n(\gamma_{i-1})\bigr| \bigr\}
\\[-2pt]
&&\qquad\quad{} + O_p \biggl(\frac{\log n}{n} \biggr)^{{(\holder+1)}/{3}}.
\end{eqnarray*}
In order to replace $B(\tau_i)$ by $B(g(\gamma_i))$, we first note that
(\ref{eqABholder}), \eqref{eqgi-fti} and Theorem~\ref{theorateUn} imply that uniformly in $i$,
\begin{eqnarray*}
\bigl|B(\tau_i)-B\bigl(g(\gamma_i)\bigr)\bigr| &\le&
C_0\bigl|\tau_i-g(\gamma_i)\bigr|^\holder
\\[-2pt]
&\le& C_0\bigl\|g'\bigr\|_{\infty}^{\holder}\bigl|f(
\tau_i)-\gamma_i\bigr|^\holder = O_p
\biggl(\frac{\log n}{n} \biggr)^{\holder/3}.
\end{eqnarray*}
By definition of $A$ and $B$, we have
$A(a)=B(g(a))|f'(g(a))|$,
for all $a\in\R$, so from Theorem~\ref{theoratefn} and~\eqref
{eqgi-Uni-1}, we conclude {that}
\begin{eqnarray*}
&& {\sup_{u\in(\tau_{i_1},\tau_{i_2}]}} B(u)\bigl\llvert \widehat{f}_n(u)-f(u)\bigr
\rrvert
\\[-2pt]
%
&&\qquad =
\max_{i_1+1\le i\le i_2} A(\gamma_i)\bigl|g(\gamma_i)-
\widehat{U}_n(\gamma_i)\bigr| \vee \max_{i_1\le i\le i_2-1} A(
\gamma_{i+1})\bigl|g(\gamma_{i+1})-\widehat{U}_n(
\gamma_i)\bigr|
\\[-2pt]
&&\qquad\quad{} + O_p \biggl(\frac{\log n}{n} \biggr)^{({\holder+1})/{3}}.
\end{eqnarray*}
By the triangle inequality, on the event $E_n$ of~\eqref{eqE} we can write
\[
|\gamma_{i+1}-\gamma_i|\le\bigl\|f'
\bigr\|_\infty \bigl\{\bigl|g(\gamma_{i+1})-\widehat {U}_n(
\gamma_{i})\bigr|+\bigl|g(\gamma_i)-\widehat{U}_n(
\gamma_i)\bigr| \bigr\}
\]
for all $i_1\le i\le i_2-1$, so that Theorem \ref{theorateUn}
together with (\ref{eqgi-Uni-1}) implies that
%
\begin{equation}
\label{eqmaxgamma} \max_{i_1\le i\le i_2-1}|\gamma_{i+1}-
\gamma_i|=O_p \biggl(\frac{\log
n}{n}
\biggr)^{1/3}.
\end{equation}
Together with (\ref{eqABholder}) and (\ref{eqgi-Uni-1}),
this allows us to replace $A(\gamma_{i+1})$ by $A(\gamma_i)$, so that
\begin{eqnarray*}
&& {\sup_{u\in(\tau_{i_1},\tau_{i_2}]}} B(u)\bigl\llvert \widehat{f}_n(u)-f(u)\bigr
\rrvert
\\[-2pt]
& &\qquad= \max_{i_1+1\le i\le i_2} A(\gamma_i)\bigl|g(\gamma_i)-
\widehat{U}_n(\gamma_i)\bigr| \vee \max_{i_1\le i\le i_2-1} A(
\gamma_{i})\bigl|g(\gamma_{i+1})-\widehat{U}_n(
\gamma_i)\bigr|
\\[-2pt]
&&\qquad\quad{} + O_p \biggl(\frac{\log n}{n} \biggr)^{({\holder+1})/{3}}.
\end{eqnarray*}
Now, recall that {$\widehat{U}_n$ is constant
on intervals $(\gamma_{i+1},\gamma_i]$, and $g$ is monotone}.
This implies that
\[
\sup_{a\in(\gamma_{i+1},\gamma_i]}\bigl|\widehat{U}_n(a)-g(a)\bigr| =\bigl |\widehat{U}_n(
\gamma_i)-g(\gamma_i)\bigr| \vee \bigl|\widehat{U}_n(
\gamma_i)-g(\gamma_{i+1})\bigr|.
\]
Therefore, taken into account joint indices, we find that
\begin{eqnarray*}
&& {\sup_{u\in(\tau_{i_1},\tau_{i_2}]}} B(u)\bigl\llvert \widehat{f}_n(u)-f(u)\bigr
\rrvert
\\[-2pt]
&&\qquad= \max_{i_1+1\le i\le i_2-1} A(\gamma_i)\sup_{a\in(\gamma_{i+1},\gamma_i]}\bigl|
\widehat{U}_n(a)-g(a)\bigr|
\\[-2pt]
&&\qquad\quad{} \vee A(\gamma_{i_2})\bigl|g(\gamma_{i_2})-\widehat{U}_n(
\gamma_{i_2})\bigr| \vee A(\gamma_{i_1})\bigl|g(\gamma_{i_1+1})-
\widehat{U}_n(\gamma_{i_1})\bigr|
\\[-2pt]
&&\qquad\quad{} + O_p \biggl(\frac{\log n}{n} \biggr)^{({\holder+1})/{3}}.
\end{eqnarray*}
Next, consider the term $A(\gamma_{i_1})|g(\gamma_{i_1+1})-\widehat
{U}_n(\gamma_{i_1})|$,
and let $\epsilon>0$.
According to~\eqref{eqcondEn} and~\eqref{eqlemtau1}, there exists
$C>0$ such that {$\mathbb{P}(I_n)>1-\epsilon$},
for $n$ sufficiently large, where
$I_n=\{\tau_{i_1}-s\le Cn^{-1/3}\}$.
By monotonicity, we have on this event that $\gamma_{i_1}=\widehat
{f}_n(\tau_{i_1})$ is between
$\widehat{f}_n(s+Cn^{-1/3})$ and~$\widehat{f}_n(s)$,
which are both equal to $f(s)+O_p(n^{-1/3})$ by~\eqref{eqEfn}.
A similar argument holds for $\gamma_{i_1+1}$, so that
%
\begin{eqnarray}
\label{eqgammai1} \gamma_{i_1}&=&f(s)+O_p
\bigl(n^{-1/3}\bigr)\quad \mbox{and}
\nonumber
\\[-9pt]
\\[-9pt]
\nonumber
 \gamma_{i_1+1}&=&f(s)+O_p
\bigl(n^{-1/3}\bigr).
\end{eqnarray}
Together with \eqref{eqgammatau} and~\eqref{eqlemEn}, this implies
\begin{eqnarray*}
\bigl|g(\gamma_{i_1+1})-\widehat{U}_n(\gamma_{i_1})\bigr| &=&
\bigl|g(\gamma_{i_1+1})-g\bigl(f(\tau_{i_1})\bigr)\bigr|
\\[-2pt]
&\le& \bigl\|g'\bigr\|_{\infty}{\bigl|\gamma_{i_1+1}-f(
\tau_{i_1})\bigr|}
\\[-2pt]
&=& \bigl\|g'\bigr\|_{\infty} {\bigl| f(s)-f(\tau_{i_1})\bigr|}+O_p
\bigl(n^{-1/3}\bigr)
\\[-2pt]
&\le& \bigl\|g'\bigr\|_{\infty} \bigl\|f'\bigr\|_{\infty}{|s-
\tau_{i_1}|}+O_p\bigl(n^{-1/3}\bigr) =
O_p\bigl(n^{-1/3}\bigr).
\end{eqnarray*}
Similarly, it follows that
%
\begin{equation}
\label{eqUn-gatgamma} \bigl|g(\gamma_{i_2})-\widehat{U}_n(
\gamma_{i_2})\bigr|=O_p\bigl(n^{-1/3}\bigr),
\end{equation}
{since by the same arguments as above,}
$\gamma_{i_2}=\hat f_n(\tau_{i_2})$ is between
$\widehat{f}_n(1-t)$ and $\widehat{f}_n(1-t-Cn^{-1/3})$
with probability {greater than $1-\epsilon$}, {and both terms} are
equal to $f(1-t)+O_p(n^{-1/3})$.
Since~$A$ is bounded, we conclude that
\begin{eqnarray*}
&& {\sup_{u\in(\tau_{i_1},\tau_{i_2}]}} B(u)\bigl\llvert \widehat{f}_n(u)-f(u)\bigr
\rrvert
\\[-2pt]
&&\qquad= \max_{i_1+1\le i\le i_2-1} A(\gamma_i)\sup_{a\in(\gamma_{i+1},\gamma_i]}\bigl|
\widehat{U}_n(a)-g(a)\bigr| \vee O_p\bigl(n^{-1/3}
\bigr)
\\[-2pt]
&&\qquad\quad{} + O_p \biggl(\frac{\log n}{n} \biggr)^{({\holder+1})/{3}}.
\end{eqnarray*}
To replace $A(\gamma_i)$ by $A(a)$ for $a\in(\gamma_{i+1},\gamma_i]$,
we use (\ref{eqmaxgamma}), (\ref{eqABholder}) and Theorem~\ref
{theorateUn}.
Together with~\eqref{eqlemi1-i2}, we conclude that
%
\begin{eqnarray}
\label{eqstep2}
\qquad && \sup_{{u\in(s,1-t]}} B(u)\bigl\llvert \widehat{f}_n(u)-f(u)
\bigr\rrvert
\nonumber
\\[-6pt]
\\[-10pt]
\nonumber
&&\qquad = \sup_{a\in(\gamma_{i_2},\gamma_{i_1{+1}}]} A(a)\bigl|\widehat{U}_n(a)-g(a)\bigr| \vee
O_p\bigl(n^{-1/3}\bigr) + O_p \biggl(
\frac{\log n}{n} \biggr)^{({\holder+1})/{3}}.
\end{eqnarray}
It remains to extend the latter supremum to the interval $[f(1-t),f(s)]$.
We have
\[
\sup_{a\in[f(1-t),\gamma_{i_2}]} A(a) \bigl|\widehat{U}_n(a)-g(a)\bigr| \le \|A
\|_{\infty} \sup_{a\in[f(1-t),\gamma_{i_2}]}\bigl|\widehat{U}_n(a)-g(a)\bigr|.
\]
According to~\eqref{eqlemmonotone},
\begin{eqnarray*}
\sup_{a\in[f(1-t),\gamma_{i_2}]}\bigl|\widehat{U}_n(a)-g(a)\bigr| &\le &\bigl|
\widehat{U}_n\bigl(f(1-t)\bigr)-g\bigl(f(1-t)\bigr)\bigr|
\\
&&{} \vee \bigl|\widehat{U}_n(\gamma_{i_2})-g(\gamma_{i_2})\bigr|
+ \bigl\|g'\bigr\|_{\infty}\bigl|\gamma_{i_2}-f(1-t)\bigr|.
\end{eqnarray*}
Similarly to~\eqref{eqgammai1}, we can write
$\gamma_{i_2}=f(1-t)+O_p(n^{-1/3})$.
Together with~\eqref{1boundUn} and~\eqref{eqUn-gatgamma} we obtain
\[
\sup_{a\in[f(1-t),\gamma_{i_2}]}\bigl|\widehat{U}_n(a)-g(a)\bigr|=O_p
\bigl(n^{-1/3} \bigr)
\]
and likewise,
\[
\sup_{a\in[\gamma_{i_1+1},f(s)]}\bigl|\widehat{U}_n(a)-g(a)\bigr|=O_p
\bigl(n^{-1/3} \bigr).
\]
From~\eqref{eqstep2}, we conclude that
%
\begin{eqnarray}
\label{eqstep3}
&&\sup_{u\in(s,1-t]} B(u)\bigl\llvert \widehat{f}_n(u)-f(u)
\bigr\rrvert
\nonumber
\\[-8pt]
\\[-8pt]
\nonumber
&&\qquad = \sup_{a\in[f(t),f(s)]} A(a)\bigl|\widehat{U}_n(a)-g(a)\bigr| \vee
R_n + O_p \biggl(\frac{\log n}{n} \biggr)^{({\holder+1})/{3}},
\end{eqnarray}
where $R_n=O_p(n^{-1/3})$.
We have
\begin{eqnarray*}
&& \P \Bigl( \sup_{a\in[f(t),f(s)]} A(a)\bigl|\widehat{U}_n(a)-g(a)\bigr| \vee
R_n \ne \sup_{a\in[f(t),f(s)]} A(a)\bigl|\widehat{U}_n(a)-g(a)\bigr|
\Bigr)
\\
&&\qquad \le \P \Bigl( R_n \ge \sup_{a\in[f(t),f(s)]} A(a)\bigr|
\widehat{U}_n(a)-g(a)\bigr| \Bigr).
\end{eqnarray*}
But it follows from Corollary \ref{corlawUn} that
%
\begin{equation}
\label{eqcorUn} \biggl( \frac{\log n}{n} \biggr)^{-1/3}
\sup_{a\in[f(t),f(s)]} A(a)\bigl|\widehat{U}_n(a)-g(a)\bigr| = 2^{-1/3} +
o_p(1).\vadjust{\goodbreak}
\end{equation}
Since $R_n=o_p((n/\log n)^{-1/3})$, it follows that the latter
probability tends to zero as $n\to\infty$.
The lemma now follows from~(\ref{eqstep3}).
\end{pf}

\begin{pf*}{Proof of Theorem \ref{theolawfn}}
Let $S_n$ be defined by~\eqref{defSn}, with $0\leq u<v\leq1$ fixed
and~$\alpha_n$ and~$\beta_n$ satisfying~\eqref{eqcondan}.
Let
%
\begin{equation}
\label{eqdefst} s_n=u+\alpha_n \quad\mbox{and}\quad
t_n=1-v+\beta_n.
\end{equation}
Then automatically $s_n$ and $t_n$ will always satisfy condition~\eqref
{eqcond73}.
If, in addition, $s_n$ and $t_n$ satisfy condition~\eqref{eqcondEn},
then according to Lemma~\ref{lemslope-inverse} together with~\eqref
{eqcorUn},
\[
\sup_{t\in(u+\alpha_n,v-\beta_n]} B(t)\bigl|\widehat{f}_n(t)-f(t)\bigr|
\]
has the same limit distribution as
\[
\sup_{a\in[f(v-\beta_n),f(u+\alpha_n)]}A(a)\bigl|\widehat{U}_n(a)-g(a)\bigr|,
\]
so that Theorem~\ref{theolawfn} follows from Corollary~\ref{corlawUn}.
When $0<u<v<1$, then $s_n$ and $t_n$ automatically satisfy~\eqref
{eqcondEn},
so we only have to consider the cases where either $u=0$ or $v=1$.
If $u=0$ and $n^{1/3}\alpha_n\to\infty$, or if $v=1$ and $n^{1/3}\beta_n\to\infty$,
then~$s_n$ and~$t_n$, as defined in~\eqref{eqdefst}, also satisfy
condition~\eqref{eqcondEn}.
Therefore, we can restrict ourselves to the case $\alpha_n=O(n^{-1/3})$
and $\beta_n=O(n^{-1/3})$.

Define $a_n=n^{-1/3}(\log n)^{1/6}$, so that $u+\alpha_n<u+a_n<v-a_n<v-\beta_n$.
By means of~\eqref{eqlemmonotone}, we find
\begin{eqnarray*}
&& \sup_{t \in(u+\alpha_n,u+a_n]} \bigl|\widehat{f}_n(t)-f(t)\bigr|
\\
&&\qquad \leq \bigl|\widehat{f}_n(u+\alpha_n)-f(u+
\alpha_n)\bigr| \vee\bigl |\widehat{f}_n(u+a_n)-f(u+a_n)\bigr|
\\
&&\qquad\quad{} + \bigl|f(u+\alpha_n)-f(u+a_n)\bigr|.
\end{eqnarray*}
By definition,
$|f(u+\alpha_n)-f(u+a_n)|\leq\|f'\|_{\infty}|\alpha_n-a_n|=O(n^{-1/3}(\log n)^{1/6})$,
and according to~\eqref{eqEfn} and~\eqref{eqEfnbord}, together
with~\eqref{eqcondan},
\begin{eqnarray*}
\widehat{f}_n(u+\alpha_n)-f(u+\alpha_n) &=&
O_p\bigl((n\alpha_n)^{-1/2}\bigr) =
o_p\bigl(n^{-1/3}(\log n)^{1/3}\bigr),
\\
\widehat{f}_n(u+a_n)-f(u+a_n) &=&
O_p\bigl(n^{-1/3}\bigr).
\end{eqnarray*}
Because $B(t)$ is uniformly bounded, it follows that
\[
\sup_{t \in(u+\alpha_n,u+a_n]} B(t)\bigl|\widehat{f}_n(t)-f(t)\bigr| = o_p
\biggl(\frac{\log n}{n} \biggr)^{1/3},
\]
and likewise
\[
\sup_{t \in(v-a_n,v-\beta_n]} B(t)\bigl|\widehat{f}_n(t)-f(t)\bigr| = o_p
\biggl(\frac{\log n}{n} \biggr)^{1/3}.
\]
This means that
\[
\sup_{t \in(u+\alpha_n,v-\beta_n]} B(t)\bigl|\widehat{f}_n(t)-f(t)\bigr| =
\sup_{t \in(u+a_n,v-a_n]} B(t)\bigl|\widehat{f}_n(t)-f(t)\bigr| \vee R_n,\vadjust{\goodbreak}
\]
where $R_n=o_p((n/\log n)^{-1/3})$.
Because $u+a_n$ and $1-v+a_n$ satisfy the conditions of Lemma~\ref
{lemslope-inverse},
together with~\eqref{eqcorUn}, it follows that
\[
\sup_{t\in(u+\alpha_n,v-\beta_n]} B(t)\bigl|\widehat{f}_n(t)-f(t)\bigr|
\]
has the same limit distribution as
\[
\sup_{c\in[f(v-a_n),f(u+a_n)]}A(c)\bigl|\widehat{U}_n(c)-g(c)\bigr|,
\]
so that Theorem~\ref{theolawfn} follows from Corollary~\ref{corlawUn}.
\end{pf*}

\section*{Acknowledgments}
The authors would like to thank Fadoua Balabdaoui for co-organizing a
research visit in Paris concerning this topic
and for helpful comments and suggestions in many stimulating discussions.

\begin{supplement}[id=suppA]
\stitle{Supplement to ``The limit distribution of the $L_\infty$-error of Grenander-type
estimators''}
\slink[doi]{10.1214/12-AOS1015SUPP} 
\sdatatype{.pdf}
\sfilename{aos1015\_supp.pdf}
\sdescription{\begin{itemize}
\item Supplement A: The supremum of the limiting process.
\item Supplement B: Preliminary results for the inverse process.
\item Supplement C: Points of jump.
\end{itemize}}
\end{supplement}

%

%

%


\printaddresses

\end{document}